\documentclass[final]{aims}
\usepackage{amsmath}
\usepackage{paralist}
\usepackage[misc]{ifsym}
\usepackage[colorlinks=true]{hyperref}
\hypersetup{urlcolor=blue,citecolor=blue}
\allowdisplaybreaks
\usepackage{url}
\def\currentvolume{X}
 \def\currentissue{X}
  \def\currentyear{200X}
   \def\currentmonth{XX}
    \def\ppages{X-XX}
     \def\DOI{10.3934/xx.xxxxxxx}
\newcommand{\MyBibitem}[2]{\bibitem{#1}}
\newcommand{\MyReference}[3]{#1. #2. #3.}
\newcommand{\MyCite}[1]{\cite{#1}}
\newcommand{\MyParagraph}[1]{\emph{#1.}\;}
\newcommand{\NN}{\mathbb{N}}
\newcommand{\ZZ}{\mathbb{Z}}
\newcommand{\RR}{\mathbb{R}}
\newcommand{\CC}{\mathbb{C}}
\newcommand{\ee}{\ensuremath{\mathrm{e}}}
\newcommand{\dd}{\ensuremath{\mathrm{d}}}
\newcommand{\ii}{\ensuremath{\mathrm{i}}}

\newcommand{\nE}{\mathcal{E}}
\newcommand{\nF}{\mathcal{F}}
\newcommand{\nL}{\mathcal{L}}
\newcommand{\nM}{\mathcal{M}}
\newcommand{\nO}{\mathcal{O}}
\newcommand{\nS}{\mathcal{S}}
\newcommand{\abs}[1]{|#1|}

\newcommand{\absbigg}[1]{\bigg|#1\bigg|}

\newcommand{\norm}[2]{\|#2\|_{#1}}
\newcommand{\normbig}[2]{\big\|#2\big\|_{#1}}
\newcommand{\nEF}{\nE^{(F)}}
\newcommand{\nSF}{\nS^{(F)}}
\newcommand{\nLF}{\nL^{(F)}}
\newcommand{\nEFh}{\nEF_h}
\newcommand{\nSFh}{\nSF_h}
\newcommand{\nLFh}{\nLF_h}
\newcommand{\nEA}{\nE^{(A)}}
\newcommand{\nEB}{\nE^{(B)}}
\newcommand{\nEAh}[1]{\nE^{(#1A)}_h}
\newcommand{\nEBh}[1]{\nE^{(#1B)}_h}
\title[Symmetric-conjugate splitting methods]{Symmetric-conjugate splitting methods for evolution equations of parabolic type} 
\author[S.~Blanes, F.~Casas, C.~Gonz{\'a}lez, M.~Thalhammer]{}
\thanks{$^*$Corresponding author: Mechthild Thalhammer}
\begin{document}
%%%%%%%%%%%%%%%%%%%%%%%%%%%%%%%%%%%%%%%%%%%%%%%%%%%%%%%%%%%%%%%%%%%%%%%%%%%%%%%%%%%%%%%%%%%%%%%%%%%%%%%%%%%%%%%%%%%%%%%%%%%%%%%%%%%%%%%%%%%%%%%%%%%
%%%%%%%%%%%%%%%%%%%%%%%%%%%%%%%%%%%%%%%%%%%%%%%%%%%%%%%%%%%%%%%%%%%%%%%%%%%%%%%%%%%%%%%%%%%%%%%%%%%%%%%%%%%%%%%%%%%%%%%%%%%%%%%%%%%%%%%%%%%%%%%%%%%
\maketitle
\centerline{\scshape
S.~Blanes$^{{\href{mailto:serblaza@imm.upv.es}{\textrm{\Letter}}}\,1}$, 
F.~Casas$^{{\href{mailto:fernando.casas@uji.es}{\textrm{\Letter}}}\,2}$, 
C.~Gonz{\'a}lez$^{{\href{mailto:cesareo.gonzalez@uva.es}{\textrm{\Letter}}}\,3}$, 
M.~Thalhammer$^{{\href{mailto:mechthild.thalhammer@uibk.ac.at}{\textrm{\Letter}}}\,*4}$
}
% Websites: \url{personales.upv.es/serblaza}, \url{www.gicas.uji.es/fernando.html}, \newline \url{www.imuva.uva.es/en/investigadores/33}, \url{techmath.uibk.ac.at/mecht}.}
\medskip
{\footnotesize
\centerline{$^1$Universitat Polit{\`e}cnica de Val{\`e}ncia, Instituto de Matem{\'a}tica Multidisciplinar, 46022~Valencia, Spain.} 
\centerline{$^2$Universitat Jaume~I, IMAC and Departament de Matem{\`a}tiques, 12071~Castell\'on, Spain.} 
\centerline{$^3$Universidad de Valladolid, Departamento de Matem{\'a}tica Aplicada, 47011~Valladolid, Spain.} 
\centerline{$^4$Universit{\"a}t Innsbruck, Institut f{\"u}r Mathematik, 6020~Innsbruck, Austria.} 
}
\bigskip
%\centerline{(Communicated by Brynjulf Owren)}
%%%%%%%%%%%%%%%%%%%%%%%%%%%%%%%%%%%%%%%%%%%%%%%%%%%%%%%%%%%%%%%%%%%%%%%%%%%%%%%%%%%%%%%%%%%%%%%%%%%%%%%%%%%%%%%%%%%%%%%%%%%%%%%%%%%%%%%%%%%%%%%%%%%
\begin{abstract}
The present work provides a comprehensive study of symmetric-conjugate operator splitting methods in the context of linear parabolic problems and demonstrates their additional benefits compared to symmetric splitting methods.
Relevant applications include nonreversible systems and ground state computations for linear Schr{\"o}dinger equations based on the imaginary time propagation.
Numerical examples confirm the favourable error behaviour of higher-order symmetric-conjugate splitting methods and illustrate the usefulness of a time stepsize control, where the local error estimation relies on the computation of the imaginary parts and thus requires negligible costs.
\end{abstract}
%%%%%%%%%%%%%%%%%%%%%%%%%%%%%%%%%%%%%%%%%%%%%%%%%%%%%%%%%%%%%%%%%%%%%%%%%%%%%%%%%%%%%%%%%%%%%%%%%%%%%%%%%%%%%%%%%%%%%%%%%%%%%%%%%%%%%%%%%%%%%%%%%%%
%%%%%%%%%%%%%%%%%%%%%%%%%%%%%%%%%%%%%%%%%%%%%%%%%%%%%%%%%%%%%%%%%%%%%%%%%%%%%%%%%%%%%%%%%%%%%%%%%%%%%%%%%%%%%%%%%%%%%%%%%%%%%%%%%%%%%%%%%%%%%%%%%%%
\section{Introduction}
%%%%%%%%%%%%%%%%%%%%%%%%%%%%%%%%%%%%%%%%%%%%%%%%%%%%%%%%%%%%%%%%%%%%%%%%%%%%%%%%%%%%%%%%%%%%%%%%%%%%%%%%%%%%%%%%%%%%%%%%%%%%%%%%%%%%%%%%%%%%%%%%%%%
%%%%%%%%%%%%%%%%%%%%%%%%%%%%%%%%%%%%%%%%%%%%%%%%%%%%%%%%%%%%%%%%%%%%%%%%%%%%%%%%%%%%%%%%%%%%%%%%%%%%%%%%%%%%%%%%%%%%%%%%%%%%%%%%%%%%%%%%%%%%%%%%%%%
%\MyParagraph{General scope}
A wide range of mathematical models for dynamical processes involves initial value problems for ordinary and partial differential equations.
Specifically designed space and time discretisation methods are of major importance in view of their effective simulation.
Over the last decades, a variety of contributions has established theoretical and numerical evidence that the class of operator splitting methods leads to favourable time integration methods and additionally preserves structural properties of linear and nonlinear evolution equations. 

\MyParagraph{General references}
The monographs~\MyCite{Kleinert1989,Messiah1999} give comprehensive overviews of applications in quantum physics.
Expositions of approved functional analytical frameworks within the scope of parabolic and Schr{\"o}dinger equations are found in~\MyCite{EngelNagel2000,HillePhillips1974,Lunardi2013,Pazy1983}.
For detailed information on splitting and composition methods, we refer to~\MyCite{HairerLubichWanner2006,McLachlanQuispel2002,SanzSernaCalvo2018}.

\MyParagraph{Objectives}
In the present work, we are concerned with the study of higher-order symmetric-conjugate splitting methods in comparison with symmetric splitting methods.
For both classes of time integration methods, certain positivity conditions on the real parts of the complex coefficients are crucial to ensure stability for nonreversible systems.
Moreover, provided that the problem data satisfy suitable regularity and consistency requirements, it is ensured that the nonstiff orders of convergence are retained.
We provide theoretical results and numerical illustrations, which confirm the reliability and efficiency of symmetric-conjugate splitting methods for relevant model problems including parabolic counterparts of linear Schr{\"o}dinger equations.
We exemplify the general observation that the numerical evolution operator associated with a symmetric-conjugate splitting method inherits the fundamental property of self-adjointness, in contrast to complex splitting methods with a strict symmetry in the configuration of their coefficients.
We conclude that a particular benefit of symmetric-conjugate splitting methods is the possibility to base local error estimations on the computation of imaginary parts.

\MyParagraph{Related works}
Our investigations are inspired by a series of works on real and complex splitting methods as well as modified splitting methods.  
As an excerpt, we mention~\MyCite{AuzingerHofstaetterKoch2019,BandraukShen1991,Bao2004,BaoDu2004,BaoJinMarkowich2002,BertoliBesseVilmart2021,BertoliVilmart2020,CaliariZuccher2021,CastellaChartierDecombesVilmart2009,Chin1997,DanailaProtas2017,DescombesDuarteDumontLaurentLouvetMassot2014,EinkemmerOstermann2016,GoldmanKaper1996,Goth2022,HansenOstermann2009a, HansenOstermann2009b,JahnkeLubich2000,Kieri2015,KozlovKvaernoOwren2004,LukassenKiehl2018,OmelyanMryglodFolk2002,OmelyanMryglodFolk2003,Sheng1989,Suzuki1991,Yoshida1990} and our former contributions~\MyCite{BaderBlanesCasas2013,BernierBlanesCasasEscorihuela2023,BlanesCasas2005,BlanesCasasMurua2010,BlanesCasasChartierMurua2013,BlanesCasasEscorihuela2022,BlanesCasasChartierEscorihuela2022,BlanesCasasGonzalezThalhammer2023a,BlanesCasasGonzalezThalhammer2023b,Thalhammer2008,Thalhammer2012}, where further references are given.

\MyParagraph{Outline}
The present work has the following structure.
In Section~\ref{sec:Section2}, we introduce operator splitting methods and give a brief summary of fundamental concepts for their convergence analysis and practical implementation.
Moreover, we introduce a unifying formulation for relevant model problems such as parabolic analogues of Schr{\"o}dinger equations.
In Section~\ref{sec:Section3}, we demonstrate the benefits of symmetric-conjugate splitting methods over symmetric splitting methods.
For this purpose, in order to reduce the amount of technicalities, it is useful to employ the generally understandable setting of real symmetric matrices.
Further explanations on the construction of higher-order splitting methods by composition and numerical comparisons of symmetric-conjugate splitting methods with standard and modified splitting methods are finally given in Section~\ref{sec:Section5}.
%%%%%%%%%%%%%%%%%%%%%%%%%%%%%%%%%%%%%%%%%%%%%%%%%%%%%%%%%%%%%%%%%%%%%%%%%%%%%%%%%%%%%%%%%%%%%%%%%%%%%%%%%%%%%%%%%%%%%%%%%%%%%%%%%%%%%%%%%%%%%%%%%%%
%%%%%%%%%%%%%%%%%%%%%%%%%%%%%%%%%%%%%%%%%%%%%%%%%%%%%%%%%%%%%%%%%%%%%%%%%%%%%%%%%%%%%%%%%%%%%%%%%%%%%%%%%%%%%%%%%%%%%%%%%%%%%%%%%%%%%%%%%%%%%%%%%%%
\section{Fundamental concepts and model problems}
\label{sec:Section2}
%%%%%%%%%%%%%%%%%%%%%%%%%%%%%%%%%%%%%%%%%%%%%%%%%%%%%%%%%%%%%%%%%%%%%%%%%%%%%%%%%%%%%%%%%%%%%%%%%%%%%%%%%%%%%%%%%%%%%%%%%%%%%%%%%%%%%%%%%%%%%%%%%%%
%%%%%%%%%%%%%%%%%%%%%%%%%%%%%%%%%%%%%%%%%%%%%%%%%%%%%%%%%%%%%%%%%%%%%%%%%%%%%%%%%%%%%%%%%%%%%%%%%%%%%%%%%%%%%%%%%%%%%%%%%%%%%%%%%%%%%%%%%%%%%%%%%%%
%\MyParagraph{Linear evolution equation}
The starting point of our considerations is the linear evolution equation 
\begin{equation}
\label{eq:IVPNonlinear}
\begin{cases}
u'(t) = F\big(u(t)\big) = A \, u(t) + B \, u(t)\,, \quad t \in [t_0, T]\,, \\
u(t_0) \text{ given}\,.
\end{cases}
\end{equation}
Throughout, we denote by $(X, \norm{X}{\cdot})$ the underlying Banach space and assume that the domains of the operators $A: D(A) \subset X \to X$ and $B: D(B) \subset X \to X$ have a non-empty intersection.
Under the presumption of suitably restricted domains, the commutator of linear operators is defined by
\begin{equation*}
\big[A, B\big] = A \, B - B A\,.
\end{equation*}  
More generally, iterated commutators are determined recursively
\begin{equation*}
\begin{gathered}
\text{ad}_A^{\, \ell}(B) = \begin{cases} B\,, & \ell = 0\,, \\ \big[A, \text{ad}_A^{\, \ell - 1}(B)\big]\,, &\ell \in \NN_{\geq 1}\,. \end{cases}
\end{gathered}
\end{equation*}

\MyParagraph{Splitting approach}
For simplicity, we restrict ourselves to uniform time grids
\begin{subequations}
\label{eq:Method}
\begin{equation}
t_n = t_0 + n h\,, \quad n \in \{0, 1, \dots, N\}\,, \quad h = \tfrac{T - t_0}{N} > 0\,, 
\end{equation}
defined by positive integer numbers $N \in \NN_{\geq 1}$.
As usual in a time-stepping approach, our aim is to compute approximations to the exact solution values through a recurrence relation of the form 
\begin{equation}
\label{eq:Recurrence}
u_{n} = \nSFh(u_{n-1}) \approx u(t_{n}) = \nEFh\big(u(t_{n-1})\big)\,, \quad n \in \{1, \dots, N\}\,, 
\end{equation}
where~$\nEF$ and~$\nSF$ represent the exact and numerical evolution operators, respectively.

Operator splitting methods rely on the idea to treat the subproblems that arise from the natural decomposition in~\eqref{eq:IVPNonlinear} separately and to compose their solutions in a favourable manner.
With regard to~\eqref{eq:Recurrence}, we consider a single subinterval and denote the evolution operators associated with 
\begin{equation*}
\begin{split}
&\begin{cases}
v'(t) = A \, u(t)\,, \quad t \in [t_{n-1}, t_n]\,, \\
v(t_{n-1}) \text{ given}\,, \quad v(t_{n}) = \nEAh{}\big(v(t_{n-1})\big)\,, 
\end{cases} \\
&\begin{cases}
w'(t) = B \, w(t)\,, \quad t \in [t_{n-1}, t_n]\,, \\
w(t_{n-1}) \text{ given}\,, \quad w(t_{n}) = \nEBh{}\big(w(t_{n-1})\big)\,, 
\end{cases}
\end{split}
\end{equation*}
by~$\nEA$ and~$\nEB$, respectively.
Incorporating suitably chosen real or more generally complex coefficients
\begin{equation}
a_j, b_j \in \CC\,, \quad j \in \{1, \dots, s\}\,,
\end{equation}
splitting methods for~\eqref{eq:IVPNonlinear} can be cast into the format~\eqref{eq:Recurrence} with 
\begin{equation}
\nSFh = \nEBh{b_s} \circ \nEAh{a_s} \circ \dots \circ \nEBh{b_1} \circ \nEAh{a_1}\,.
\end{equation}
Throughout, we denote by $p \in \NN$ the classical order of a splitting method and tacitly assume that the coefficients of the considered schemes fulfill the elementary consistency condition
\begin{equation}
\label{eq:Consistency}
\sum_{j=1}^{s} a_j = 1\,, \quad \sum_{j=1}^{s} b_j = 1\,.
\end{equation}
On occasion, when the structural characteristics are essential, we use the compact symbolic notation 
\begin{equation}
\nSFh = h \, \big(b_s, a_s, \dots, b_1, a_1\big)\,.  
\end{equation}
\end{subequations}

\MyParagraph{Symmetric-conjugate methods}
We primarily focus on symmetric-conjugate operator splitting methods of the form 
\begin{subequations}
\label{eq:SCSplitting}
\begin{equation}
\begin{gathered}
s = 2 \, r\,, \quad a_1 = 0\,, \\
a_{s+2-j} = \overline{a_j}\,, \quad j \in \{2, \dots, r\}\,, \quad b_{s+1-j} = \overline{b_j}\,, \quad j \in \{1, \dots, r\}\,, \\
\nSFh = h \, \big(\overline{b_1}, \overline{a_2}, \overline{b_2}, \dots, \overline{a_r}, \overline{b_r}, a_{r+1}, b_r, a_r, \dots, b_2, a_2, b_1, 0\big)\,,
\end{gathered}
\end{equation}
and impose the positivity conditions
\begin{equation}
\begin{gathered}
a_j \in \CC\,, \quad \Re(a_j) > 0\,, \quad j \in \{2, \dots, r+1\}\,, \\
b_j \in \CC\,, \quad \Re(b_j) > 0\,, \quad j \in \{1, \dots, r\}\,, 
\end{gathered}
\end{equation}
\end{subequations}
to ensure well-definedness and thus stability for evolution equations of parabolic type.
It is natural to contrast symmetric-conjugate with symmetric splitting methods 
\begin{equation}
\label{eq:SSplitting}
\begin{gathered}
s = 2 \, r\,, \quad a_1 = 0\,, \\
a_{s+2-j} = a_j\,, \quad j \in \{2, \dots, r\}\,, \quad b_{s+1-j} = b_j\,, \quad j \in \{1, \dots, r\}\,, \\
\nSFh = h \, \big(b_1, a_2, b_2, \dots, a_r, b_r, a_{r+1}, b_r, a_r, \dots, b_2, a_2, b_1, 0\big)\,.
\end{gathered}
\end{equation}
For notational simplicity, we henceforth omit coefficients that are equal to zero and write for instance 
\begin{equation*}
\begin{gathered}
\nSFh = h \, \big(\overline{b_1}, \overline{a_2}, \overline{b_2}, \dots, \overline{a_r}, \overline{b_r}, a_{r+1}, b_r, a_r, \dots, b_2, a_2, b_1\big)\,, \\
\nSFh = h \, \big(b_1, a_2, b_2, \dots, a_r, b_r, a_{r+1}, b_r, a_r, \dots, b_2, a_2, b_1\big)\,,
\end{gathered}
\end{equation*}
for short, see~\eqref{eq:SCSplitting} and~\eqref{eq:SSplitting}.
The different behaviour exhibited by both classes of splitting methods and the particularly favourable performance of symmetric-conjugate schemes in the time integration of linear ordinary differential equations that are defined by real symmetric matrices deserves a detailed analysis, which is carried out in this work.

\MyParagraph{Elementary splitting methods}
As elementary instances, we introduce the famous Lie--Trotter and Strang splitting methods and a third-order symmetric-conjugate splitting method.   
We recall that the positive integers $s \in \NN_{\geq 1}$ and $p \in \NN_{\geq 1}$ denote the number of stages and the classical order of a splitting method~\eqref{eq:Method}. 
\begin{enumerate}[(i)]
\item 
The simplest first-order scheme
\begin{equation*}
\begin{gathered}
p = 1\,, \quad s = 1\,, \quad a_1 = 1\,, \quad b_1 = 1\,, \\
\nSFh = \nEBh{} \circ \nEAh{} \approx \nEFh\,,
\end{gathered}
\end{equation*}
is known as Lie--Trotter splitting method.
Evidently, it fulfils the positivity condition $a_1, b_1 > 0$, but it does not fit into the classes~\eqref{eq:SCSplitting} or~\eqref{eq:SSplitting}, respectively. 
\item
The second-order Strang splitting method, which comprises two stages and the symmetric composition 
\begin{equation}
\label{eq:Strang}
\begin{gathered}
p = 2\,, \quad s = 2\,, \quad a_1 = 0\,, \quad a_2 = 1\,, \quad b_1 = \tfrac{1}{2} = b_2\,, \\
\nSFh = \nEBh{\frac{1}{2}} \circ \nEAh{} \circ \nEBh{\frac{1}{2}} \approx \nEFh\,, 
\end{gathered}
\end{equation}
is contained in both classes~\eqref{eq:SCSplitting} and~\eqref{eq:SSplitting} with $r = 1$.
\item
The simplest symmetric-conjugate splitting method of order three 
\begin{equation*}
\begin{gathered}
p = 3\,, \quad s = 3\,, \quad a_1 = 0\,, \quad a_2 = \tfrac{1}{2} \, \big(1 + \ii  \, \tfrac{1}{\sqrt{3}}\big)\,, \quad b_1 = \tfrac{1}{2} \, a_2\,, \quad b_2 = \tfrac{1}{2}\,, \\
\nSFh = \nEBh{\overline{b_1}} \circ \nEAh{\overline{a_2}} \circ \nEBh{b_2} \circ \nEAh{a_2} \circ \nEBh{b_1} \approx \nEFh\,, 
\end{gathered}
\end{equation*}
was proposed in~\MyCite{BandraukShen1991}.
Alternatively, it is retained as a special double jump composition of the Strang splitting method
\begin{equation*}
\nSFh = \nS^{(F, \, \text{Strang})}_{\overline{a_2} \, h} \circ \nS^{(F, \, \text{Strang})}_{a_2 \, h}
= h \, \big(\tfrac{1}{2} \, \overline{a_2}, \overline{a_2}, \tfrac{1}{2} \, (a_2 + \overline{a_2}), a_2, \tfrac{1}{2} \, a_2\big)\,.
\end{equation*}
\end{enumerate}

\MyParagraph{Higher-order splitting methods}
In our numerical tests, detailed in Section~\ref{sec:Section5}, we compare higher-order standard and modified operator splitting methods involving real coefficients with complex symmetric and symmetric-conjugate splitting methods, see~Figure~\ref{fig:Figure1}. 
For the convenience of the readers, we display the method coefficients of the symmetric-conjugate schemes with corresponding denominations in Figures~\ref{fig:Figure2} and~\ref{fig:Figure3}.
Besides, a link to a Matlab code is provided in Section~\ref{sec:Section5}.

A fourth-order symmetric scheme by \text{Yoshida}~\MyCite{Yoshida1990} comprises negative coefficients, wherefore instabilities arise for evolution equations of parabolic type, see~\eqref{eq:SchemesOrder4RealSymmetric}.
Its complex analogue given in~\eqref{eq:SchemesOrder4ComplexSymmetric} leads to a stable alternative.
A non-standard scheme is \textsc{Chin}’s fourth-order modified potential operator splitting method~\MyCite{Chin1997} involving positive coefficients and double commutators, see also~\MyCite{BlanesCasasGonzalezThalhammer2023a,BlanesCasasGonzalezThalhammer2023b}.

Moreover, we apply different complex symmetric and symmetric-conjugate splitting methods, which satisfy the stability condition that all coefficients have non-negative real parts.
Optimised symmetric schemes of order four and a symmetric scheme of order six are found in~\MyCite{BlanesCasasChartierMurua2013}.\footnote{See also~\url{www.gicas.uji.es/Research/splitting-complex.html}.}
Amongst the symmetric-conjugate splitting methods, we highlight sixth-order schemes, recently proposed in~\MyCite{BernierBlanesCasasEscorihuela2023} for linear unitary problems.
They set aside the second- and fourth-order barriers for standard and modified splitting methods, see~\MyCite{AuzingerHofstaetterKoch2019,BlanesCasas2005,Sheng1989,Suzuki1991,GoldmanKaper1996} and references given therein.
We point out that these schemes are characterised by positive coefficients $a_j > 0$ for $j \in \{2, \dots, s\}$.
Consequently, they are suitable for the time integration of different classes of evolution equations including parabolic and Schr{\"o}dinger equations, see~\eqref{eq:ParabolicLinear} as well as~\eqref{eq:SchroedingerLinear} below and Table~\ref{tab:Table1}.

\MyParagraph{Model problems}
As prototype models for linear evolution equations of parabolic type, we study partial differential equations that involve the Laplacian and a potential.
This kind of nonreversible systems in particular arises in ground and excited state computations based on the imaginary time propagation, see~\MyCite{Kleinert1989,Messiah1999} for detailed information on the theoretical foundations.

In the context of these relevant applications, it is often reasonable to assume that solutions are localised and highly regular. 
Thus, it is justified to restrict the originally unbounded space domain to a Cartesian product of sufficiently large intervals and to impose periodic boundary conditions.
Conveniences of such settings are that fast Fourier techniques can be applied and that higher-order splitting methods do not suffer from severe order reductions.
Throughout, we base our considerations on these assumptions.

In the following, we denote by $\Omega \subseteq \RR^d$ the underlying space domain, by $\Delta = \partial_{x_1}^2 + \dots + \partial_{x_d}^2$ the Laplacian with respect to the spatial variables $x = (x_1, \dots, x_d) \in \Omega$, and by $V: \Omega \to \RR$ a space-dependent real-valued potential acting as multiplication operator. 
For notational simplicity, we omit scaling constants and signs, unless they are significant with regard to classifications as parabolic or Schr{\"o}dinger equations.

As a test problem, we consider the linear parabolic problem 
\begin{equation}
\label{eq:ParabolicLinear}  
\begin{cases}
\partial_t U(x, t) = \tfrac{1}{2} \, \Delta \, U(x, t) - V(x) \, U(x, t)\,, \\
U(x, t_0) \text{ given}\,, \quad (x, t) \in \Omega \times [t_0, T]\,, 
\end{cases}
\end{equation}  
for a real-valued solution $U: \Omega \times [t_0, T] \to \RR$. 
This special choice is justified by the imaginary time propagation of the linear Schr{\"o}dinger equation
\begin{equation}
\label{eq:SchroedingerLinear}  
\begin{cases}
\ii \, \partial_t \Psi(x, t) = - \, \tfrac{1}{2} \, \Delta \Psi(x, t) + V(x) \, \Psi(x, t)\,, \\
\Psi(x, t_0) \text{ given}\,, \quad (x, t) \in \Omega \times [t_0, T]\,, 
\end{cases}
\end{equation}  
with complex-valued wave function $\Psi: \Omega \times [t_0, T] \to \CC$, see Section~\ref{sec:Section5} for further explanations.

Evidently, the above model problems~\eqref{eq:ParabolicLinear} and~\eqref{eq:SchroedingerLinear} can be cast into the unifying formulation
\begin{subequations}
\label{eq:ModelProblems}  
\begin{equation}
\label{eq:PDE}  
\begin{cases}
&\partial_t U(x, t) = \alpha \, \Delta \, U(x, t) + \beta \, V(x) \, U(x, t)\,, \\
&U(x, t_0) \text{ given}\,, \quad (x, t) \in \Omega \times [t_0, T]\,, 
\end{cases}
\end{equation}
with $\alpha\,, \beta \in \CC$ denoting certain constants and $U: \Omega \times [t_0, T] \to \CC$ the solution.
Specifically, we chose the arising quantities as
\begin{gather}
\label{eq:LinearP}  
\alpha = \tfrac{1}{2}\,, \quad \beta = - \, 1\,, \\ 
\label{eq:LinearS}  
\alpha = \tfrac{1}{2} \, \ii\,, \quad \beta = - \, \ii\,. 
\end{gather}
\end{subequations}
Additional scaling constants and signs that are insignificant with regard to classifications as parabolic or Schr{\"o}dinger equations are again neglected.

\MyParagraph{Practicable presumptions}
In this work, we are primarily interested in the time integration of parabolic initial-boundary value problems by high-order operator splitting methods.
Hence, to ensure that the nonstiff orders of convergence are retained in our numerical comparisons, we restrict ourselves to situations, where the problem data satisfy suitable regularity and consistency requirements.
Otherwise, even though well-tailored higher-order schemes are generally more favourable than lower-order schemes, we have to expect substantial order reductions.

Furthermore, to implement the physically relevant and numerically challenging case of three space dimensions by means of Fourier spectral space discretisations, we simplify the general setting in this respect. 
Specifically, we presume that it is appropriate to replace the underlying space domain with a Cartesian product of sufficiently large intervals
\begin{subequations}
\label{eq:Fourier}  
\begin{equation}
\label{eq:Omega}  
\Omega = [- \, a, a]^d\,, \quad a > 0\,, \quad d \in \{1, 2, 3\}\,, 
\end{equation}  
and prescribe Gaussian-like initial states
\begin{equation}
U(x, t_0) = c_1 \, \ee^{\, - \, c_2 \, \abs{x - c_3}^2}\,, \quad x = (x_1, \dots, x_d) \in \Omega\,, 
\end{equation}  
which fulfil intrinsic periodicity conditions on~$\Omega$ with high accuracy.

In the following, we sketch the employed means that permit efficient implementations of operator splitting methods for the model problem~\eqref{eq:PDE} by fast Fourier techniques.
For detailed descriptions, we refer the readers to our former works \MyCite{BlanesCasasGonzalezThalhammer2023a,BlanesCasasGonzalezThalhammer2023b}.

\MyParagraph{Fourier series representations}
We denote by $(\nF_m)_{m \in \ZZ^d}$ the Fourier functions with periodicity domain~\eqref{eq:Omega} and by $(\lambda_m)_{m \in \ZZ^d}$ the corresponding real eigenvalues of the Laplace operator 
\begin{equation}
\begin{gathered}
\nF_m(x) = (2 \, a)^{- \, \frac{d}{2}} \; \ee^{\, \ii \, \pi \, m_1 \, (\frac{x_1}{a} + 1)} \cdots \, \ee^{\, \ii \, \pi \, m_d \, (\frac{x_d}{a} + 1)}\,, \\
\Delta \, \nF_m = \lambda_m \, \nF_m\,, \quad \lambda_m = - \, \tfrac{\pi^2 \abs{m}^2}{a^2} \in \RR\,, \\
x = (x_1, \dots, x_d) \in \Omega\,, \quad m = (m_1, \dots, m_d) \in \ZZ^d\,.
\end{gathered}  
\end{equation}  
Realisations of representations by Fourier series 
\begin{equation}
\begin{gathered}
v = \sum_{m \in \ZZ^d} v_m \, \nF_m\,, \quad \Delta v = \sum_{m \in \ZZ^d} \lambda_m \, v_m \, \nF_m\,, \\
v_m = \int_{[- \, a, a]^d} v(x) \, \nF_{-m}(x) \; \dd x\,, \quad m \in \ZZ^d\,, 
\end{gathered}  
\end{equation}  
\end{subequations}
are based on suitable truncations of the infinite index sets $\nM \subset \ZZ^d$ such that $\abs{\nM} = M \in \NN$ as well as quadrature approximations by the trapezoidal rule. 

\MyParagraph{Stiff and nonstiff subproblems}
We rewrite the partial differential equation in~\eqref{eq:PDE} as abstract evolution equation of the form~\eqref{eq:IVPNonlinear} with $u(t) = U(\cdot, t)$ for $t \in [t_0, T]$.
We employ the natural decomposition into a stiff and a nonstiff partial differential equation
\begin{equation*}
\begin{cases}
&\partial_t U(x, t) = \alpha \, \Delta \, U(x, t)\,, \\
&\partial_t U(x, t) = \beta \, V(x) \, U(x, t)\,.
\end{cases}
\end{equation*}
Here, it should be noted that the considered Fourier spectral space discretisation affects the definition of the operators.

Accordingly, we assign the first unbounded linear operator with the Laplacian
\begin{subequations}
\label{eq:SubproblemASolution}  
\begin{equation}
A = \alpha \, \Delta\,, \quad \alpha \in \RR\,, \quad \alpha > 0\,.
\end{equation}
For any complex coefficient $a \in \CC$ with non-negative real part $\Re(a) \geq 0$, it is ensured that the corresponding subproblem 
\begin{equation}
\begin{cases}
v'(t) = a \, A \, v(t)\,, \quad t \in [t_{n-1}, t_n]\,, \\
\displaystyle v(t_{n-1}) = \sum_{m \in \ZZ^d} v_m(t_{n-1}) \, \nF_m \text{ given}\,, 
\end{cases} 
\end{equation}
is well-posed and its solution formally given by the Fourier series representation 
\begin{equation}
v(t_n) = \nEAh{a} \, v(t_{n-1}) = \sum_{m \in \ZZ^d} \ee^{\, a \, h \, \alpha \, \lambda_m} \, v_m(t_{n-1}) \, \nF_m\,,
\end{equation}  
\end{subequations}  
see also~\eqref{eq:Fourier}.

The second nonstiff subproblem is defined by the potential, which acts as a multiplication operator, and resolved by pointwise products.

\MyParagraph{Evolution equations of Schr{\"o}dinger type}
It is instructive to observe that linear Schr{\"o}dinger equations~\eqref{eq:SchroedingerLinear} are included in~\eqref{eq:PDE}, see also~\eqref{eq:LinearS}.
For this reason, we may expect that the convergence analysis provided in~\MyCite{Thalhammer2008,Thalhammer2012} for Fourier spectral space discretisations combined with high-order operator splitting methods based on real coefficients can be transferred to evolution equations of parabolic type and complex symmetric-conjugate splitting methods satisfying a positivity condition.
For the convenience of the readers, we next recall fundamental notions and concepts.

\MyParagraph{Local and global errors}
We employ a standard argument based on the telescopic identity to conclude that the validity of stability bounds combined with local error expansions implies global error estimates of the form  
\begin{equation*}
\normbig{X}{u_n - u(t_n)} \leq C \, \big(\normbig{X}{u_0 - u(t_0)} + h^p\big)\,, \quad n \in \{1, \dots, N\}\,.
\end{equation*}
This fundamental principle serves as guideline for our convergence analysis of operator splitting methods applied to evolution equations of parabolic type. 

The general approach is most comprehensible within the context of evolution equations that are defined by bounded linear operators 
\begin{equation*}
u'(t) = F\big(u(t)\big) = (A + B) \, u(t)\,, \quad t \in [t_0, T]\,, \quad A, B: X \longrightarrow X\,,
\end{equation*}
since then the exact and splitting solutions can be represented by products of matrix exponentials
\begin{equation*}  
\nEFh = \ee^{\, h \, (A + B)}\,, \quad
\nSFh = \ee^{\, b_s h B} \, \ee^{\, a_s h A} \cdots \ee^{\, b_1 h B} \, \ee^{\, a_1 h A}\,, 
\end{equation*}
and thus stability bounds follow at once from 
\begin{equation*}  
\begin{gathered}
\normbig{X \leftarrow X}{\nEFh} \leq \ee^{\, C h}\,, \quad \normbig{X \leftarrow X}{\nSFh} \leq \ee^{\, C h}\,, \\
C = \max\Bigg\{\norm{X \leftarrow X}{A} + \norm{X \leftarrow X}{B}, \sum_{j=1}^{s} \Big(\abs{a_j} \norm{X \leftarrow X}{A} + \abs{b_j} \norm{X \leftarrow X}{B}\Big)\bigg\}\,.
\end{gathered}
\end{equation*}
Provided that the coefficients satisfy certain order conditions such that the local error expansion 
\begin{equation*}
\nLFh = \nSFh - \nEFh = \nO\big(h^{p+1}\big)
\end{equation*}
holds, the desired global error estimate follows by means of the relations 
\begin{equation*}
\begin{gathered}
u_n - u(t_n) = \big(\nSFh\big)^n \, \big(u_0 - u(t_0)\big) 
+ \sum_{k=0}^{n-1} \big(\nSFh\big)^{n-1-k} \nLFh \big(\nEFh\big)^k \, u(t_0)\,, \\
\normbig{X}{u_n - u(t_n)} \leq \ee^{\, C t_n} \, \Big(\normbig{X}{u_0 - u(t_0)} + n \, \normbig{X \leftarrow X}{\nLFh} \normbig{X}{u(t_0)}\Big)\,, \\
n \in \{1, \dots, N\}\,.
\end{gathered}
\end{equation*}
%%%%%%%%%%%%%%%%%%%%%%%%%%%%%%%%%%%%%%%%%%%%%%%%%%%%%%%%%%%%%%%%%%%%%%%%%%%%%%%%%%%%%%%%%%%%%%%%%%%%%%%%%%%%%%%%%%%%%%%%%%%%%%%%%%%%%%%%%%%%%%%%%%%
%%%%%%%%%%%%%%%%%%%%%%%%%%%%%%%%%%%%%%%%%%%%%%%%%%%%%%%%%%%%%%%%%%%%%%%%%%%%%%%%%%%%%%%%%%%%%%%%%%%%%%%%%%%%%%%%%%%%%%%%%%%%%%%%%%%%%%%%%%%%%%%%%%%
\section{Symmetric-conjugate versus symmetric methods}
\label{sec:Section3}
%%%%%%%%%%%%%%%%%%%%%%%%%%%%%%%%%%%%%%%%%%%%%%%%%%%%%%%%%%%%%%%%%%%%%%%%%%%%%%%%%%%%%%%%%%%%%%%%%%%%%%%%%%%%%%%%%%%%%%%%%%%%%%%%%%%%%%%%%%%%%%%%%%%
%%%%%%%%%%%%%%%%%%%%%%%%%%%%%%%%%%%%%%%%%%%%%%%%%%%%%%%%%%%%%%%%%%%%%%%%%%%%%%%%%%%%%%%%%%%%%%%%%%%%%%%%%%%%%%%%%%%%%%%%%%%%%%%%%%%%%%%%%%%%%%%%%%%
%\MyParagraph{Structural properties of complex splitting methods}
In this section, we contrast the favourable properties of symmetric-conjugate splitting methods comprising complex coefficients with those of symmetric splitting methods.
Essential ingredients are series expansions that characterise local errors and the spectral theorem. 
The construction of higher-order schemes by composition and numerical illustrations are described in Section~\ref{sec:Section5}.

\MyParagraph{Restrictions and generalisations}
We point out that the following analysis based on infinite series expansions is powerful regarding the treatment of high-order splitting methods, but it also has some restrictions. 

Our main conclusions concerning the accumulation of inaccurate imaginary parts rely on the assumption that the defining operators correspond to real symmetric matrices.

To a certain extent, our setting is associated with spatial semi-discretisations of partial differential equations, but it should be noted that the constants $C_A, C_B$ in~\eqref{eq:CACB} below increase when the space grids are refined.
Furthermore, for parabolic evolution equations, the inverses of~$\nE$ and~$\nS$ in~\eqref{eq:ES1} below are not well-defined, since they involve evaluations at negative times, e.g. 
\begin{equation*}
\big(\nE(t)\big)^{-1} = \nE(- \, t)\,, \quad t \in (0, T - t_0]\,.
\end{equation*}
We indicate this issue for the simplest splitting method, the Lie--Trotter splitting method.
In situations, where~$A$ represents the Laplacian and~$B$ a multiplication operator related to a real-valued potential, the numerical evolution operator and its time derivative
\begin{equation*}
\nS(t) = \ee^{\, t A} \, \ee^{\, t B}\,, \quad \tfrac{\dd}{\dd t} \, \nS(t) = \ee^{\, t A} \, (A + B) \, \ee^{\, t B}\,, \quad t \in (0, T - t_0]\,, 
\end{equation*}
are well-defined.  
However, the inverse
\begin{equation*}
\big(\nS(t)\big)^{-1} = \ee^{- \, t B} \, \ee^{- \, t A}\,, \quad t \in (0, T - t_0]\,,
\end{equation*}
and likewise the composition  
\begin{equation*}
G(t) = \tfrac{\dd}{\dd t} \, \nS(t) \, \big(\nS(t)\big)^{-1} = \ee^{\, t A} \, (A + B) \, \ee^{- \, t A}\,, \quad t \in (0, T - t_0]\,, 
\end{equation*}
involve negative times and hence necessitate the restriction to finite dimensional cases, which result from spatial semi-discretisations.  

A rigorous analysis in the lines of~\MyCite{Thalhammer2008,Thalhammer2012} is intended for future work.
It will be based on suitable adaptations of the arguments and will require alternative solution representations as well as specifications of the employed expansions and arising remainders.
Fundamental means for stepwise expansions of the exact and numerical evolution operators will be provided by the variation-of-constants formula and Taylor series expansions.
The characterisation of the resulting regularity requirements will be linked to the identification of iterated commutators.

\MyParagraph{Simplified setting of matrices}
We consider the initial value problem for a linear ordinary differential equation
\begin{subequations}
\label{eq:CACB}  
\begin{equation}
\begin{cases}
u'(t) = F\big(u(t)\big) = A \, u(t) + B \, u(t)\,, \quad t \in [t_0, T]\,, \\
u(t_0) \in \RR^M \text{ given}\,,
\end{cases} 
\end{equation}
under the additional assumption that the defining matrices are real and symmetric 
\begin{equation}
\begin{gathered}
A \in \RR^{M \times M}\,, \quad \norm{}{A} \leq C_A\,, \quad A^{*} = A^T = A\,, \\
B \in \RR^{M \times M}\,, \quad \norm{}{B} \leq C_B\,, \quad B^{*} = B^T = B\,.
\end{gathered}
\end{equation}
\end{subequations}
The exact and numerical evolution operators
\begin{equation}  
\label{eq:ES1}
\begin{gathered}
\nE_t = \nE(t) = \ee^{\, t \, (A + B)} \in \RR^{M \times M}\,, \\ 
\nS_t = \nS(t) = \ee^{\, b_s \, t \, B} \, \ee^{\, a_s \, t \, A} \cdots \, \ee^{\, b_1 \, t \, B} \, \ee^{\, a_1 \, t \, A} \in \CC^{M \times M}\,, \\
t \in [0, T - t_0]\,,
\end{gathered}
\end{equation}
are given by exponential series, that is 
\begin{equation*}  
\ee^{\, t \, L} = \sum_{\ell=0}^{\infty} \tfrac{1}{\ell!} \, t^{\ell} L^{\ell}\,, \quad \normbig{}{\ee^{\, t \, L}} \leq \ee^{\, \abs{t} \, \norm{}{L}}\,, \quad t \in \RR\,.
\end{equation*}
We note that a splitting method involving complex coefficients yields complex-valued approximations to the real-valued solution and recall the consistency condition~\eqref{eq:Consistency} ensuring order $p \in \NN_{\geq 1}$.
For the ease of notation, we only indicate the dependence on the time increment.

\MyParagraph{Series representations}
In the context of matrices, we may use formal infinite series expansions of decisive components to deduce substantive results for splitting methods with complex coefficients.

On the one hand, the exact evolution operator
\begin{equation*}
\nE(t) = \ee^{\, t \, (A + B)}\,, \quad t \in [0, T - t_0]\,, 
\end{equation*}
satisfies the initial value problem 
\begin{equation*}
\begin{cases} \nE'(t) = (A + B) \, \nE(t)\,, \quad t \in [0, T - t_0]\,, \\ \nE(0) = I\,, \end{cases}
\end{equation*}
where $I \in \RR^{M \times M}$ denotes the identity matrix. 

On the other hand, we make use of the fact that the evolution operator associated with a splitting method is given as the exponential of a time-dependent operator 
\begin{equation}
\label{eq:OmegaH}
\nS(t) = \ee^{\, \Omega(t)} = \ee^{\, t \, H(t)}\,, \quad t \in [0, T - t_0]\,, 
\end{equation}
and fulfills a related nonautonomous linear differential equation 
\begin{equation*}
\begin{cases} \nS'(t) = G(t) \, \nS(t)\,, \quad t \in [0, T - t_0]\,, \\ \nS(0) = I\,, \end{cases}
\end{equation*}
where formally $G = \nS' \, \nS^{-1}$ as well as $\Omega' = (\exp(\text{ad}_{\Omega}) - I)^{-1} \, \text{ad}_{\Omega} \, G$.
Important findings are that statements on the difference $H - (A + B)$ allow to draw conclusions on $G - (A + B)$ and hence on~$\nS - \nE$.

More precisely, formal representations of $H - (A + B)$ as infinite series can be found by applying recursively the Baker--Campbell--Hausdorff formula to the numerical evolution operator~\eqref{eq:ES1}.
Considering the Lie algebra $\nL(A,B)$ generated by $\{A,B\}$ with the commutator
\begin{equation*}
[A, B] = A \, B - B A
\end{equation*}
as Lie bracket and denoting by $\nL_{\ell}(A, B)$ the homogeneous subspace of degree $\ell \in \NN_{\geq 2}$ with $k$th basis element $E_{\ell k}(A, B)$ for $k \in K_{\ell}$, e.g.
\begin{equation*}
\begin{gathered}
E_{21}(A, B) = [A,B]\,, \quad K_2 = \{1\}\,, \\
E_{31}(A, B) = \big[A, [A, B]\big]\,, \quad E_{32}(A, B) = \big[B, [A, B]\big]\,, \quad K_3 = \{1, 2\}\,.
\end{gathered}
\end{equation*}
we obtain the formal series expansion 
\begin{equation}
\label{eq:H}
H(t) = A + B + \sum_{\ell = 2}^{\infty} t^{\ell - 1} \sum_{k \in K_{\ell}} e_{\ell k}(a, b) \, E_{\ell k}(A, B)\,, \quad t \in [0, T - t_0]\,,
\end{equation}
where $e_{\ell k}(a, b)$ represents a polynomial of degree~$k$ with respect to the complex coefficients $(a_j, b_j)_{j=1}^{s}$.

\MyParagraph{Particular structures of symmetric-conjugate schemes}
The above formal representation~\eqref{eq:H} is valid for arbitrary real matrices $A, B \in \RR^{M \times M}$ and splitting methods with complex coefficients. 
Provided that~$A$ and~$B$ are symmetric and the considered splitting methods are symmetric-conjugate, it turns out that~$H$ has a particular structure. 

We henceforth fix the time increment $h = \tfrac{T - t_0}{N}$ for some positive integer $N \in \NN_{\geq 1}$ and use again the convenient notation $\nE_h = \nE(h)$, $\nS_h = \nS(h)$, and $H_h = H(h)$.
For brevity, we do not indicate the dependencies on the defining operators $E_{\ell k} = E_{\ell k}(A, B)$ and the splitting coefficients $e_{\ell k} = e_{\ell k}(a, b)$.

Evidently, the properties real and symmetric are inherited by the exact evolution operator
\begin{equation*}  
\nE_h = \nE_h^T \in \RR^{M \times M}\,.
\end{equation*}
Moreover, for symmetric-conjugate splitting methods, we conclude that
\begin{equation}  
\label{eq:HSelfAdjoint}
\nS_h = \nS_h^{*} \in \CC^{M \times M}\,, \quad H_h = H_h^{*} \in \CC^{M \times M}\,,
\end{equation}
see also~\eqref{eq:SCSplitting} and~\eqref{eq:CACB}. 
Note, however, that this is in general not true for symmetric splitting methods. 
Observing that the iterated commutators of real symmetric matrices satisfy the relations 
\begin{equation*}
E_{\ell k}^T = (-1)^{\ell + 1} \, E_{\ell k}\,, \quad k \in K_{\ell}\,, \quad \ell \in \NN_{\geq 2}\,,
\end{equation*}
e.g. $E_{21}^T = [A , B]^T = B A - A B = - \, E_{21}$, the difference of the formal series expansion~\eqref{eq:H} and its adjoint   
\begin{equation*}
\begin{split}
0 &= H_h - H_h^{*} = \sum_{\ell = 2}^{\infty} h^{\ell - 1} \sum_{k \in K_{\ell}} \big(e_{\ell k} \, E_{\ell k} - \overline{e_{\ell k}} \, E_{\ell k}^T\big) \\
&= \sum_{\ell = 2}^{\infty} h^{\ell - 1} \sum_{k \in K_{\ell}} \Big(\big(1 + (-1)^{\ell}\big) \, \Re(e_{\ell k}) + \ii \, \big(1 + (-1)^{\ell + 1}\big) \, \Im(e_{\ell k})\Big) \, E_{\ell k}\,,
\end{split}
\end{equation*}
imply that the arising coefficients are either real or purely imaginary quantities 
\begin{equation*}
\begin{cases} e_{\ell k} = r_{\ell k}\,, & \text{$\ell$ odd}\,, \\ e_{\ell k} = \ii \, r_{\ell k}\,, & \text{$\ell$ even}\,, \end{cases} \quad r_{\ell k} \in \RR\,, \quad k \in K_{\ell}\,, \quad \ell \in \NN_{\geq 2}\,.
\end{equation*}
In consequence, we obtain a decomposition of the complex matrix~$H_h$ into real symmetric and skew-symmetric contributions
\begin{equation}
\label{eq:HSymmetric}
\begin{gathered}
H_h = A + B + H_h^{(\RR, \text{sym})} + \ii \, H_h^{(\RR, \text{skew})}\,, \\
H_h^{(\RR, \text{sym})} = \underset {\ell \text{ odd}}{\sum_{\ell = 3}^{\infty}} h^{\ell - 1} \sum_{k \in K_{\ell}} r_{\ell k} \, E_{\ell k} \in \RR^{M \times M}\,, \\
H_h^{(\RR, \text{sym})} = \big(H_h^{(\RR, \text{sym})}\big)^T\,, \\
H_h^{(\RR, \text{skew})} = \underset{\ell \text{ even}}{\sum_{\ell = 2}^{\infty}} h^{\ell - 1} \sum_{k \in K_{\ell}} r_{\ell k} \, E_{\ell k} \in \RR^{M \times M}\,, \\
H_h^{(\RR, \text{skew})} = - \, \big(H_h^{(\RR, \text{skew})}\big)^T\,.
\end{gathered}
\end{equation}

\MyParagraph{Errors in imaginary parts}
The subsequent considerations will explain the favourable behaviour of symmetric-conjugate splitting methods in comparison with symmetric splitting methods for linear ordinary differential equations that are defined by real symmetric matrices.
We first restate the above decomposition~\eqref{eq:HSymmetric}
 in terms of local errors and next analyse the accumulation of errors over time. 
\begin{enumerate}[(i)]
\item     
A splitting method has order $p \in \NN_{\geq 1}$, i.e.~$\nS_h - \nE = \nO(h^{p+1})$, if $H_h - (A + B) = \nO(h^p)$. 
Specifically, if a symmetric-conjugate operator is of even order $p \in \NN_{\geq 2}$, then
\begin{equation*}
H_h^{(\RR, \text{sym})} = \nO\big(h^p\big)\,, \quad H_h^{(\RR, \text{skew})} = \nO\big(h^{p+1}\big)\,.
\end{equation*}
Otherwise, if it is of odd order $p \in \NN_{\geq 1}$, then
\begin{equation*}
H_h^{(\RR, \text{sym})} = \nO\big(h^{p+1}\big)\,, \quad H_h^{(\RR, \text{skew})} = \nO\big(h^p\big)\,.
\end{equation*}
\item
By means of the spectral theorem applied to the self-adjoint matrix $H_h = H_h^{*} \in \CC^{M \times M}$, see~\eqref{eq:HSelfAdjoint}, we conclude that there exist a unitary matrix $U_h \in \CC^{M \times M} $ and a real diagonal matrix $D_h \in \RR^{M \times M}$ such that
\begin{equation*}
H_h = U_h \, D_h \, U_h^{*}\,. 
\end{equation*}
Notice that evaluation at $h = 0$ implies that the corresponding transformation matrix is real 
\begin{equation*}
A + B = H_0 = U_0 \, D_0 \, U_0^T\,, \quad U_0 \in \RR^{M \times M}\,. 
\end{equation*}
Due to $H_h = H_0 + \nO(h^p)$, the matrices depend smoothly on the time increment
\begin{equation*}
U_h = U_0 + \nO\big(h^p\big)\,, \quad D_h = D_0 + \nO\big(h^p\big)\,.
\end{equation*}
In particular, the errors in the imaginary parts fulfill
\begin{equation*}
\Im(U_h) = \nO\big(h^p\big)\,.
\end{equation*}
In consequence, we obtain the following identities for the exact and numerical evolution operators and multiple compositions 
\begin{equation}
\label{eq:ES}  
\begin{gathered}
\nE_h = \ee^{\, h \, (A + B)} = U_0 \, \ee^{\, h \, D_0} \, U_0^T\,, \\
\nE_h^n = \big(U_0 \, \ee^{\, h \, D_0} \, U_0^T\big)^{n} = U_0 \, \ee^{\, n \, h \, D_0} \, U_0^T\,, \\
\nS_h = \ee^{\, h \, H_h} = U_h \, \ee^{\, h \, D_h} \, U_h^{*}\,, \\
\nS_h^n = \big(U_h \, \ee^{\, h \, D_h} \, U_h^{*}\big)^{n} = U_h \, \ee^{\, n \, h \, D_h} \, U_h^{*}\,, \\
\quad n \in \{0, 1, \dots, N\}\,. 
\end{gathered}
\end{equation}
As $\ee^{\, n \, h \, D_h}$ is bounded by $\ee^{\, T \, \norm{}{D_h}}$, we finally conclude that the relation 
\begin{equation}
\label{eq:dyn1}
\nS_h^n = U_0 \, \ee^{\, n \, h \, D_h} \, U_0^T + \nO\big(h^p\big)\,, \quad n \in \{0, 1, \dots, N\}\,, 
\end{equation}
is valid, where the implicit constant in the $\mathcal{O}$ term does not depend on $n$.
\end{enumerate}

\MyParagraph{Conclusions}
The identity in~\eqref{eq:dyn1} has several remarkable consequences, which we confirm and complement by numerical illustrations for model problems with real-valued solutions in Section~\ref{sec:Section5}. 
\begin{enumerate}[(i)]
\item
When applied to the initial value problem
\begin{equation}
\label{eq:ODEAB}
\begin{cases} u'(t) = (A + B) \, u(t)\,, \quad t \in [t_0, T]\,, \\ u(t_0) = u_0\,, \end{cases}
\end{equation}
with time increment $h = \tfrac{T - t_0}{N}$, any symmetric-conjugate splitting method of order $p \in \NN_{\geq 1}$ is conjugate to the exact solution to the initial value problem
\begin{equation*}
\begin{cases} v'(t) = D_h \, v(t)\,, \quad t \in [t_0, T]\,, \\ v(t_0) = u_0\,, \end{cases}
\end{equation*}
where the real diagonal matrix~$D_h = D_0 + \nO(h^p)$ is a perturbation of the same order of the matrix diagonalising the real symmetric matrix $A + B$, see~\eqref{eq:SCSplitting} and~\eqref{eq:CACB}.
\item
The numerical approximation to the real-valued solution has an imaginary part of the same order~$p$, when $p$ is odd, or of order $p+1$, when $p$ is even, respectively, see also Table~\ref{tab:Table2}.
\item
This error does not accumulate, and, hence, it does not affect the global performance. 
More precisely, for symmetric-conjugate splitting methods, the relative errors in the imaginary part
\begin{equation*}
\frac{\norm{}{\Im(u_n)}}{\norm{}{u_n}}\,, \quad u_n = \nS_h^n u_0\,, \quad n \in \{0, 1, \dots, N\}\,,
\end{equation*}
remain bounded over time, since they are only due to the transformation by $U_h = U_0 + \nO(h^p)$.
For a comparison of fourth-order symmetric-conjugate versus symmetric splitting methods, we refer to Figure~\ref{fig:Model1}.
The considered situation is related to a ground state computation by the imaginary time propagation for a linear Schrödinger equation under a quartic potential.
We mirror the errors in the imaginary parts and the ground state energy.
\end{enumerate}
Altogether, we conclude that symmetric-conjugate splitting methods are particularly favourable for the numerical approximation of linear ordinary differential equations that are defined by real symmetric matrices~\eqref{eq:CACB}.
The incorporation of complex coefficients requires the use of complex arithmetics and increases the computational effort, but it permits the design of high-order schemes, that satisfy the stability condition~\eqref{eq:SCSplitting} and thus overcome the second-order and fourth-order barriers for standard and modified splitting methods. 
Illustrative numerical results are given in Section~\ref{sec:Section5}.
%%%%%%%%%%%%%%%%%%%%%%%%%%%%%%%%%%%%%%%%%%%%%%%%%%%%%%%%%%%%%%%%%%%%%%%%%%%%%%%%%%%%%%%%%%%%%%%%%%%%%%%%%%%%%%%%%%%%%%%%%%%%%%%%%%%%%%%%%%%%%%%%%%%
%%%%%%%%%%%%%%%%%%%%%%%%%%%%%%%%%%%%%%%%%%%%%%%%%%%%%%%%%%%%%%%%%%%%%%%%%%%%%%%%%%%%%%%%%%%%%%%%%%%%%%%%%%%%%%%%%%%%%%%%%%%%%%%%%%%%%%%%%%%%%%%%%%%
\section{Test equations and numerical comparisons}
\label{sec:Section5}
%%%%%%%%%%%%%%%%%%%%%%%%%%%%%%%%%%%%%%%%%%%%%%%%%%%%%%%%%%%%%%%%%%%%%%%%%%%%%%%%%%%%%%%%%%%%%%%%%%%%%%%%%%%%%%%%%%%%%%%%%%%%%%%%%%%%%%%%%%%%%%%%%%%
%%%%%%%%%%%%%%%%%%%%%%%%%%%%%%%%%%%%%%%%%%%%%%%%%%%%%%%%%%%%%%%%%%%%%%%%%%%%%%%%%%%%%%%%%%%%%%%%%%%%%%%%%%%%%%%%%%%%%%%%%%%%%%%%%%%%%%%%%%%%%%%%%%%
%\MyParagraph{Implementation}
In this section, we provide numerical evidence confirming and complementing our theoretical analysis of symmetric-conjugate splitting methods for spatial semi-discretisations of parabolic equations.

A Matlab code, which illustrates the practical implementation of operator splitting methods combined with Fourier spectral space discretisations for three-dimensional model problems is available at
\begin{center}
\url{doi.org/10.5281/zenodo.8238819}.
\end{center}
It in particular reproduces numerical results presented in the sequel.

In connection with our first illustration and the numerical computation of global errors, we focus on the case of a single space dimension.
However, for schemes comprising negative coefficients, severe stability issues have to be expected when the spatial grid width is refined or in higher dimensions, since then the problems become significantly stiffer.

\MyParagraph{Symmetric-conjugate versus symmetric methods}
A first numerical illustration is related to the relevant issue of ground state computations for quantum-mechanical systems. 
The provided comparisons for different splitting methods of order four verify and complement our observations in Section~\ref{sec:Section3}.
In particular, they show that symmetric-conjugate splitting methods possess distinctive features for evolution equations that are defined by self-adjoint operators and have real-valued solutions. 
Contrary to symmetric counterparts, it is ensured that the numerical evolution operators inherit the property of self-adjointness, which results in favourable approximations and is also reflected in the errors in the imaginary parts and the ground state energies. 

\emph{Test equation.} 
We consider the linear Schr{\"o}dinger equation~\eqref{eq:SchroedingerLinear} for a quartic potential and a localised Gaussian-like initial state.  
Moreover, we impose periodic boundary conditions on a sufficiently large space interval such that the effect of the truncation error is insignificant.
A related parabolic problem is obtained by integration in imaginary time, i.e.~by replacing the time variable with $- \, \ii \, t$, see also~\eqref{eq:ParabolicLinear} and~\eqref{eq:PDE}.
For the convenience of the readers, we restate the test equation and specify the data of the problem 
\begin{equation}
\label{eq:TestProblem1}  
\begin{gathered}
\begin{cases}
&\partial_t U(x, t) = \tfrac{1}{2} \, \Delta \, U(x, t) - V(x) \, U(x, t)\,, \\
&U(x, t_0) = U_0(x)\,, \quad (x, t) \in \Omega \times [0, T]\,.
\end{cases} \\
V: \RR \longrightarrow \RR: x \longmapsto 5 - \tfrac{1}{2} \, x^2 + \tfrac{1}{80} \, x^4\,, \\
U_0: \RR \longrightarrow \RR: x \longmapsto \tfrac{1}{\sqrt[4]{\pi}} \, \ee^{- \frac{1}{2} (x-1)^2}\,, \\
\Omega = [- \, a, a]\,, \quad a = 10\,, \quad M = 256\,.
\end{gathered}
\end{equation}
Within our setting, the eigenvalues of the discrete Laplace operator are non-positive, see~\eqref{eq:Fourier}.
Moreover, due to the fact that the potential takes non-negative values, it is ensured that the discretisation of the defining operator $\frac{1}{2} \, \Delta - V$ has negative eigenvalues.
The time propagation of~\eqref{eq:TestProblem1} combined with suitable projection thus yields approximations to stationary states of the quantum-mechanical system.
The ground state is linked to the lowest energy level in modulus, and excited states correspond to higher energies in modulus. 

Basically, the spatially discretised system is defined by real and symmetric matrices $A, B \in \RR^{M \times M}$, see~\eqref{eq:CACB} and~\eqref{eq:ODEAB}, respectively, and hence the following considerations are appropriate.
The eigenvalues of $A + B$ are negative numbers and the solution value at the final time can formally be written as a linear combination of associated normalised eigenvectors 
\begin{equation*}
\begin{gathered}
(A + B) \, v_m = E_m \, v_m\,, \\
E_m \in \RR_{< 0}\,, \quad v_m \in \RR^M\,, \quad \norm{}{v_m} = 1\,, \quad m \in \{0, 1, \dots, M-1\}\,, \\
u(T) = \sum_{m=0}^{M-1} c_m \, \ee^{\, T \, E_m} \, v_m\,.
\end{gathered}
\end{equation*}
Provided that the dominant eigenvalue is simple, i.e.~$E_m < E_0$ for $m \in \{1, \dots, M-1\}$, and the corresponding coefficient nonzero $c_0 \neq 0$ such that
\begin{equation*}
u(T) = c_0 \, \ee^{\, T \, E_0} \bigg(v_0 + \sum_{m=1}^{M-1} \tfrac{c_m}{c_0} \, \ee^{- \, T \, (E_0 - E_m)} \, v_m\bigg) \approx c_0 \, \ee^{\, T \, E_0} \, v_0\,,
\end{equation*}
this allows to determine a numerical approximation to the first eigenvector related to the ground state 
\begin{equation*}
\tfrac{1}{\norm{}{u(T)}} \, u(T) \approx v_0\,.
\end{equation*}
We point out that an appropriate choice of the final time has to be adjusted to the location of the eigenvalues, i.e.~the smaller the difference $E_0 - E_1$, the larger $T > 0$.
The computation of the ground state energy then relies on the identity 
\begin{equation*}
E_0 = v_0^T (A + B) \, v_0\,.
\end{equation*}
The space discretisation and in particular the features of the time discretisation method will affect the quality of the obtained numerical results.

\begin{figure}[t!]
\begin{center}
\includegraphics[width=0.7\textwidth]{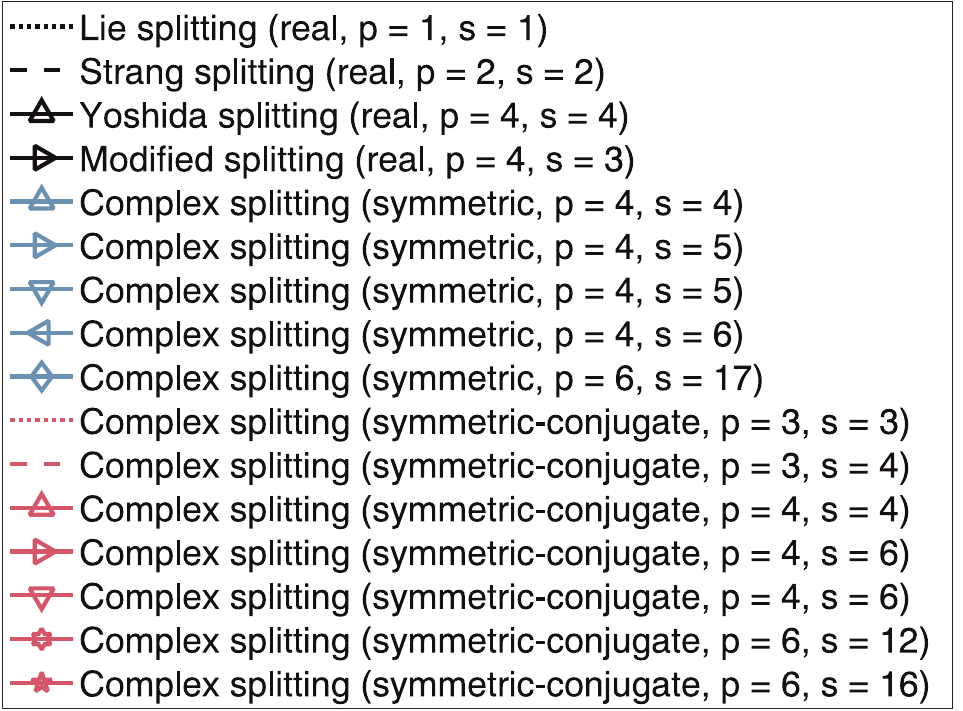}
\caption{Real and complex splitting methods applied in numerical tests. Denominations and characteristics (nonstiff order~$p$, number of stages $s$).
The coefficients of the symmetric-conjugate schemes are given in Figures~\ref{fig:Figure2} and~\ref{fig:Figure3}.}
\label{fig:Figure1}
\end{center}
\end{figure}

\begin{figure}[t!]
\begin{center}
\includegraphics[width=0.95\textwidth]{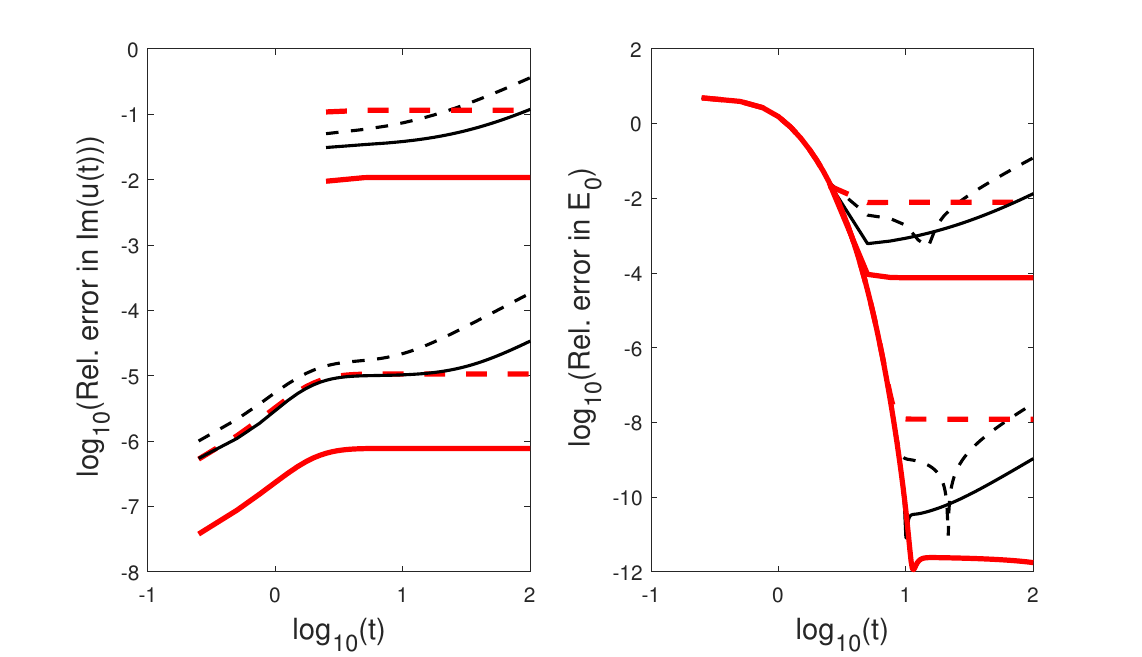} 
\caption{Time integration of the parabolic model problem~\eqref{eq:TestProblem1} by non-optimised and optimised fourth-order operator splitting methods involving complex coefficients with time increments $h = \frac{T}{40}$ (top curves) and $h=\frac{T}{400}$ (bottom curves).  
Relative errors in the imaginary parts of the numerical solutions over time (left) and corresponding errors in the ground state energy (right).
Symmetric schemes comprising $s = 4$ (thin black dashed line) and $s > 4$ (thin black solid line) stages with increasing errors in a log-log scale
versus symmetric-conjugate schemes comprising $s = 4$ (thick red dashed line) and $s > 4$ (thick red solid line) stages with bounded errors.}
\label{fig:Model1}
\end{center}
\end{figure}

\emph{Fourth-order splitting methods.} 
For the time integration of~\eqref{eq:TestProblem1}, we apply symmetric and symmetric-conjugate splitting methods of order four, respectively, see Figures~\ref{fig:Figure1} and~\ref{fig:Figure2}.

For the purpose of illustration, we comment on the construction of the fourth-order schemes comprising four stages from the second-order Strang splitting method
\begin{equation*}
\nSFh = \nS^{[2, \RR]}_h = h \, \big(\tfrac{1}{2}, 1, \tfrac{1}{2}\big)\,, 
\end{equation*}
by means of the composition technique, see also~\eqref{eq:Method} and~\eqref{eq:Strang}.
In order to distinguish different schemes, we adapt our former notation and indicate the numbers of stages as well as structural characteristics, but omit the defining function.
Specifically, we use the triple jump composition 
\begin{equation*}
\begin{split}  
\nS^{[4]}_h &= \nS_{\alpha_3 \, h}^{[2]} \circ \nS_{\alpha_2 \, h}^{[2]} \circ  \nS_{\alpha_1 \, h}^{[2]} \\
&= h \, \big(\tfrac{1}{2} \, \alpha_3, \alpha_3, \tfrac{1}{2} \, (\alpha_2 + \alpha_3), \alpha_2, \tfrac{1}{2} \, (\alpha_2 + \alpha_1), \alpha_1, \tfrac{1}{2} \, \alpha_1\big)\,.
\end{split}
\end{equation*}
Provided that the arising coefficients fulfil the conditions 
\begin{equation}
\label{eq:alphas}
\sum_{j=1}^3 \alpha_j = 1\,, \quad \sum_{j=1}^3 \alpha_j^3 = 0\,,
\end{equation}
this leads to splitting methods of nonstiff order four. 
These algebraic equations have the real-valued solutions
\begin{subequations}
\label{eq:SchemesOrder4}
\begin{equation}
\label{eq:SchemesOrder4RealSymmetric}
\begin{gathered}
\gamma = 2 - 2^{\frac{1}{3}}\,, \quad \alpha_1 = \gamma^{-1}\,, \quad \alpha_2 = 1 - 2 \, \alpha_1\,, \quad \alpha_3 = \alpha_1\,, \\
\nS^{[4, \RR, \text{sym}]}_h = h \, \big(\tfrac{\alpha_1}{2}, \alpha_1, \tfrac{1 - \alpha_1}{2}, 1 - 2 \, \alpha_1, \tfrac{1-\alpha_1}{2}, \alpha_1, \tfrac{\alpha_1}{2}\big)\,, 
\end{gathered}
\end{equation}
which define the standard symmetric \text{Yoshida} splitting method, see~\MyCite{Yoshida1990}.
In addition, they admit a complex-valued solution, which correspond to the symmetric splitting method 
\begin{equation}
\label{eq:SchemesOrder4ComplexSymmetric}
\begin{gathered}
\gamma = 2 - 2^{\frac{1}{3}} \, \ee^{\, \frac{2}{3} \, \ii \, \pi}\,, \quad
\alpha_1 = \gamma^{-1}\,, \quad \alpha_2 = 1 - 2 \, \alpha_1\,, \quad \alpha_3 = \alpha_1\,, \\
\nS^{[4, \CC, \text{sym}]}_h = h \, \big(\tfrac{\alpha_1}{2}, \alpha_1, \tfrac{1 - \alpha_1}{2}, 1 - 2 \, \alpha_1, \tfrac{1-\alpha_1}{2}, \alpha_1, \tfrac{\alpha_1}{2}\big)\,, 
\end{gathered}
\end{equation}
as well as a complex-valued solution, which correspond to the symmetric-conjugate splitting method
\begin{equation}
\begin{gathered}
\alpha_1 = \tfrac{1}{4} + \ii \, \tfrac{\sqrt{15}}{12}\,, \quad \alpha_2 = \tfrac{1}{2}\,, \quad \alpha_3 = \overline{\alpha_1}\,, \\
\nS^{[4, \CC, \text{sym-conj}]}_h
= h \, \big(\tfrac{\overline{\alpha_1}}{2}, \overline{\alpha_1}, \tfrac{\overline{\alpha_1} + \alpha_2}{2}, \alpha_2, \tfrac{\alpha_1 + \alpha_2}{2}, \alpha_1, \tfrac{\alpha_1}{2}\big)\,.
\end{gathered}
\end{equation}
\end{subequations}
Further schemes are obtained from~\eqref{eq:alphas} by complex conjugation.
The symmetric and symmetric-conjugate splitting methods have in common that they contain the same number of exponentials and provide fourth-order approximations. 
Moreover, the sizes of the main error terms at order five, measured as
\begin{equation*}
\text{err} = \absbigg{\sum_{j=1}^3 \alpha_j^5}\,,  
\end{equation*}
are about 200 times smaller than the error of the triple jump composition with real coefficients
\begin{equation*}
\text{err}^{[4, \RR, \text{sym}]}_h \approx 5.3\,, \quad
\text{err}^{[4, \CC, \text{sym}]}_h \approx 0.024\,, \quad
\text{err}^{[4, \CC, \text{sym-conj}]}_h \approx 0.028\,, 
\end{equation*}
see also~\MyCite{BlanesCasasEscorihuela2022}.

The additionally considered symmetric and symmetric-conjugate splitting methods of order four comprise a higher number of stages $s > 4$ and include positive coefficients $(a_j)_{j=1}^{s}$, which distinguishes them from the schemes with four stages and ensures stability for parabolic problems as well as Schr{\"o}dinger equations.  
The selected optimised symmetric-conjugate splitting method shows a favourable accuracy behaviour and leads to relatively small errors.

\emph{Numerical results.} 
In order to confirm the expected qualitative and quantitative differences between the above described fourth-order symmetric-conjugate and symmetric splitting methods, we prescribe a certain time increment $h > 0$ and perform the integration until a sufficiently large final time $T = N h$ is reached.
In addition, a reference solution $u(T) \in \RR^{M \times M}$ with real values up to machine precision and the associated ground state energy $E_0 = u(T)^T (A + B) \, u(T)$ are computed numerically.
Then, at each time step, the relative sizes of the imaginary parts with respect to the solution values and the relative errors of the ground state energies 
\begin{equation*}
\begin{gathered}
\frac{\|\Im(u_n)\|}{\|u_n\|}\,, \\
\frac{\abs{E_0 - \Re(u_n)^T (A + B) \, \Re(u_n)}}{\abs{E_0}}\,, \\
n \in \{0, 1, \dots, N\}\,, 
\end{gathered}
\end{equation*}
are determined.

As explained in Section~\ref{sec:Section3}, for the symmetric-conjugate splitting methods, it is expected that the errors in the imaginary parts remain bounded, whereas for symmetric splitting methods involving complex coefficients a significant error growth over time may occur.
Concerning the errors in the approximation of the ground state energy, the quantities $\Re(u_n)^T(A + B) \, \Re(u_n)$ for $n \in \{0, 1, \dots, N\}$ constitute approximations to the ground state of a perturbed matrix that depends on the splitting method and the time increment. 
Hence, it is expected that the errors decrease up to a certain time, which depends on the order of the method and the time increment.
Beyond that time, it is again assumed that the error accumulation in the imaginary part will lead to a significant error growth for symmetric splitting methods and bounded errors for symmetric-conjugate schemes. 

The numerical results displayed in Figure~\ref{fig:Model1} illustrate how the different characters of symmetric-conjugate versus symmetric splitting methods manifest in practice.
In the left panel, we depict the relative errors in the imaginary parts of the numerical solutions over time. 
In the right panel, we display the corresponding relative errors in the ground state energy. 
We indeed observe that the symmetric schemes introduce errors that grow linearly in a log-log scale, which may lead to unphysical effects. 
Contrary, after a transient time, the errors committed by the symmetric-conjugate schemes are nearly constant.

\MyParagraph{Local and global errors}
\begin{subequations}
In the lines of our first illustration, we next consider the linear test equation~\eqref{eq:PDE} subject to periodic boundary conditions on a suitably chosen space interval and apply the operator splitting methods listed in Figure~\ref{fig:Figure1}.
Specifically, we set $\alpha = 1 = \beta$ and study the quadratic and quartic potentials
\begin{equation}
V: \RR \longrightarrow \RR: x \longmapsto x^2\,, \quad V: \RR \longrightarrow \RR: x \longmapsto \tfrac{1}{24} \, x^4\,. 
\end{equation}
The space grid points and the time interval are given by
\begin{equation}
\begin{gathered}
\Omega = [- \, a, a]\,, \quad a = 10\,, \quad M = 100\,, \quad t_0 = 0\,, \quad T \in \{1, 10\}\,.
\end{gathered}
\end{equation}
For the problem, when considered on the whole real line with quadratic potential and initial state
\begin{equation}
U_0: \RR \longrightarrow \RR: x \longmapsto \ee^{- \frac{1}{2} \, x^2}\,,
\end{equation}
the exact solution is known.
In this special situation, we thus have the possibility to verify the correctness of the implementation and conclude that the errors caused by the truncation of the space domain, the implicitely imposed periodic boundary conditions, and the application of the Fourier spectral discretisation method are insignificant.
In the general case, we determine the local and global time discretisation errors with respect to numerical reference solutions.
\end{subequations}

\begin{figure}[t!]
\begin{center}
\includegraphics[width=0.49\textwidth]{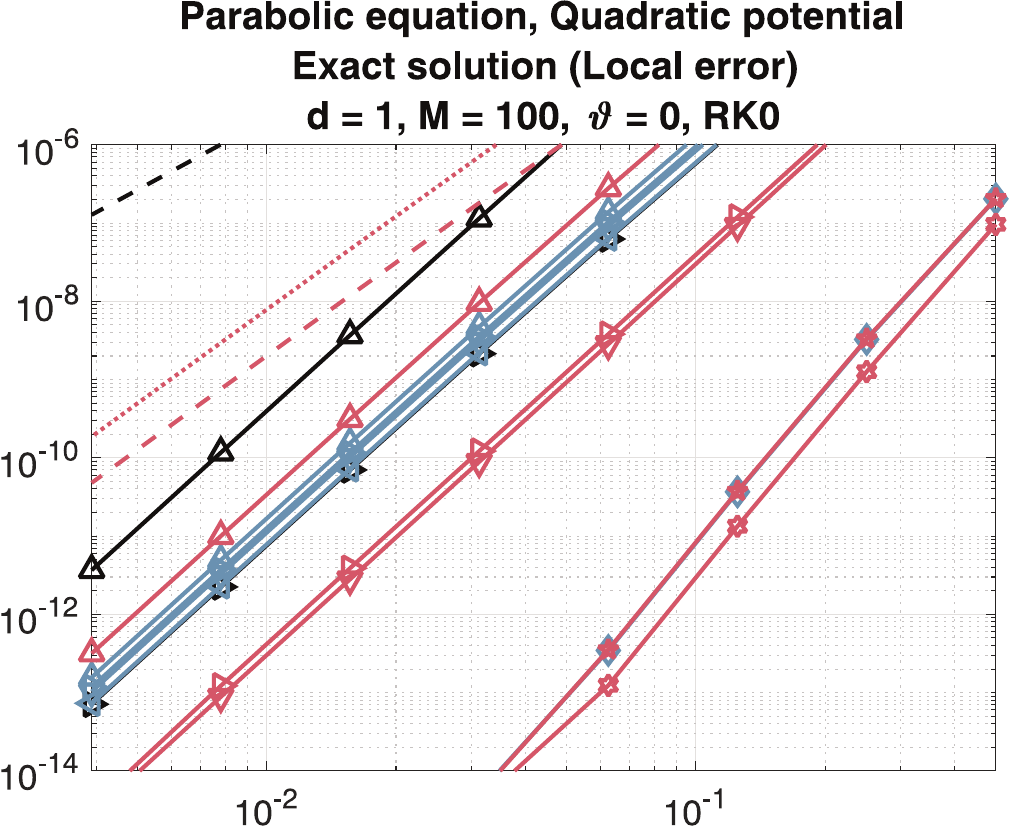} 
\includegraphics[width=0.49\textwidth]{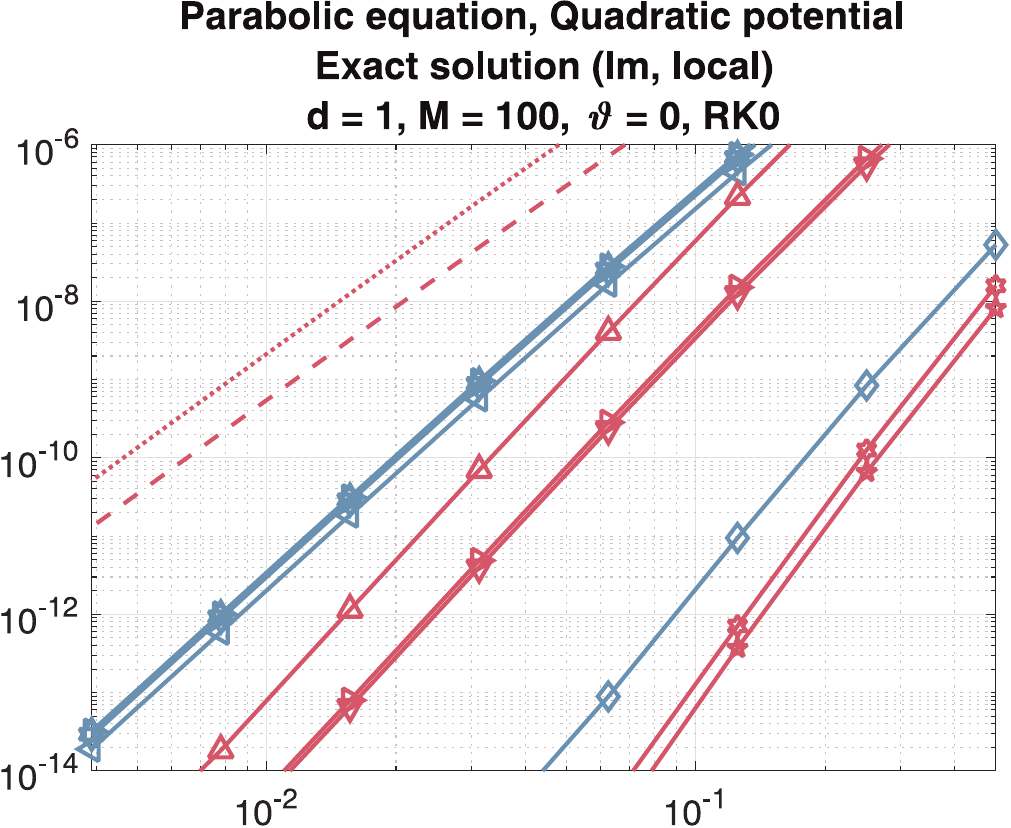} \\[2mm]
\includegraphics[width=0.49\textwidth]{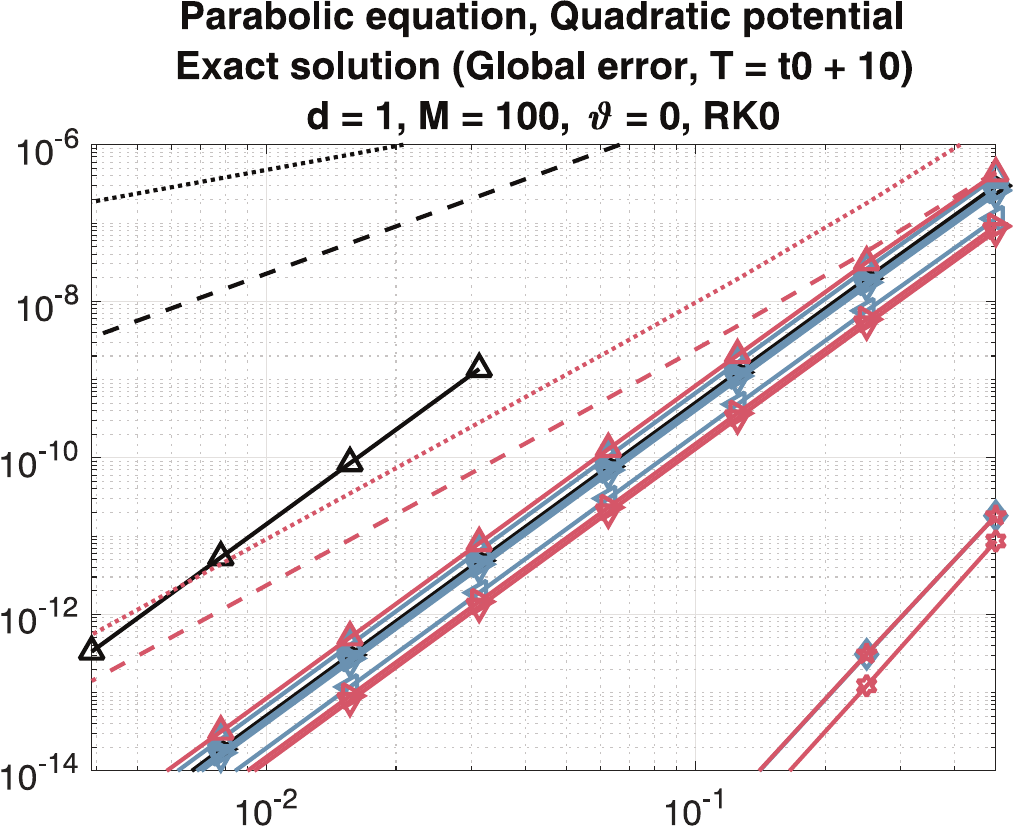} 
\includegraphics[width=0.49\textwidth]{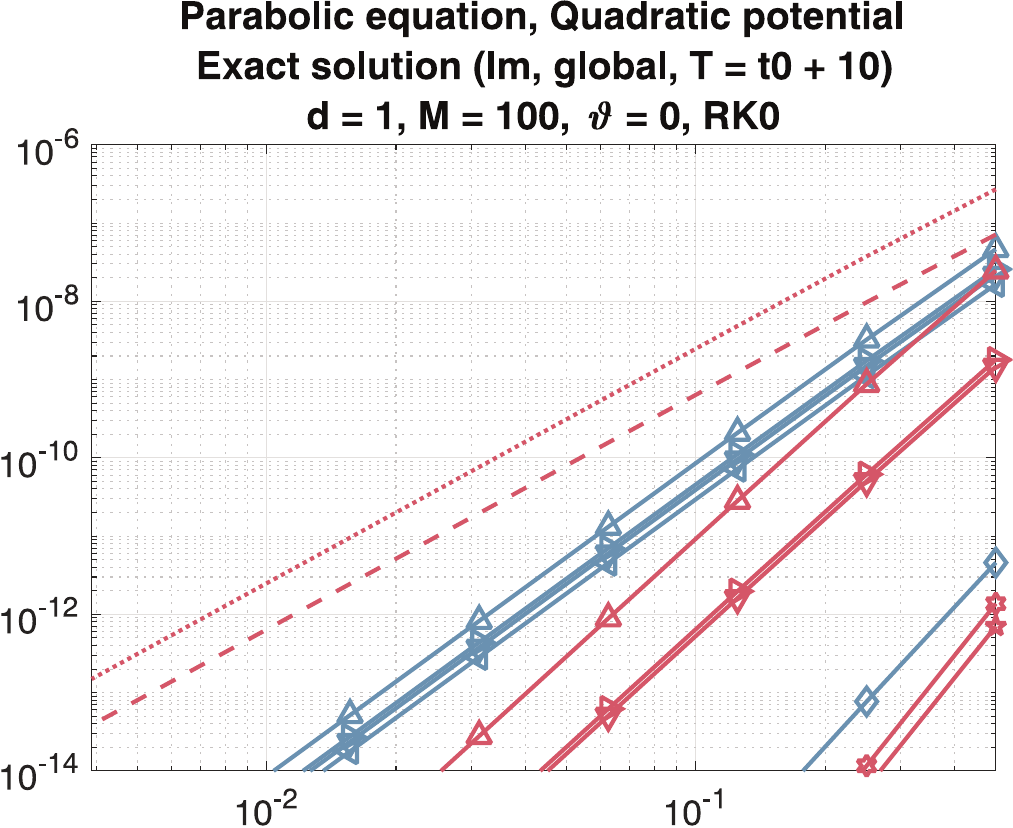} 
\caption{Time integration of the linear parabolic model problem with real-valued solution by real and complex splitting methods, see also~\eqref{eq:ModelProblems} and Figure~\ref{fig:Figure1}.
For the considered quadratic potential, the exact solution is known.
Left: Local and global errors.
Right: Corresponding errors in the imaginary parts.}
\label{fig:Model2_1}
\end{center}
\end{figure}

\begin{figure}[t!]
\begin{center}
\includegraphics[width=0.49\textwidth]{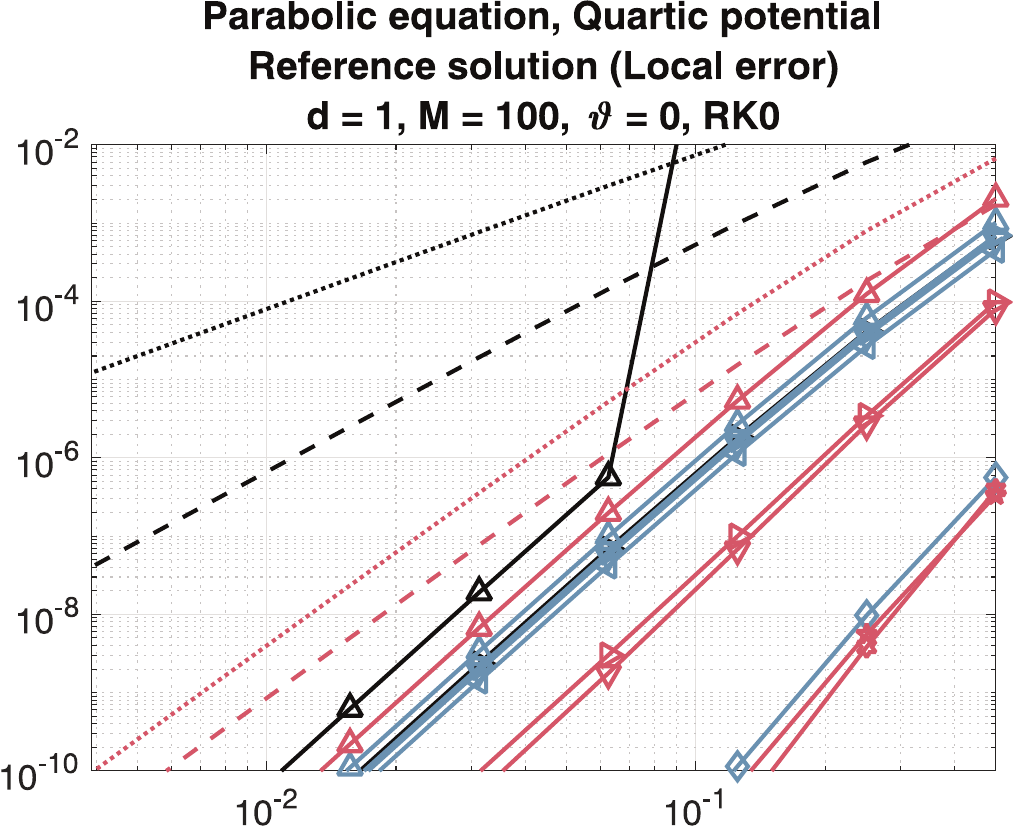} 
\includegraphics[width=0.49\textwidth]{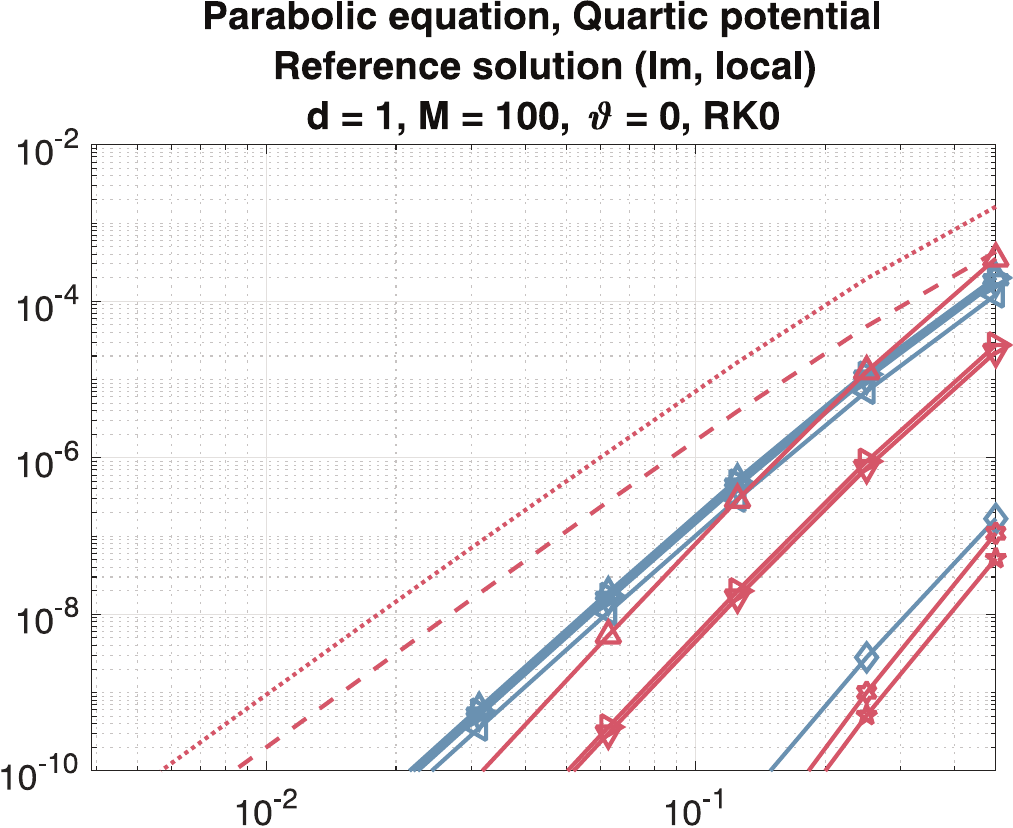} \\[2mm]
\includegraphics[width=0.49\textwidth]{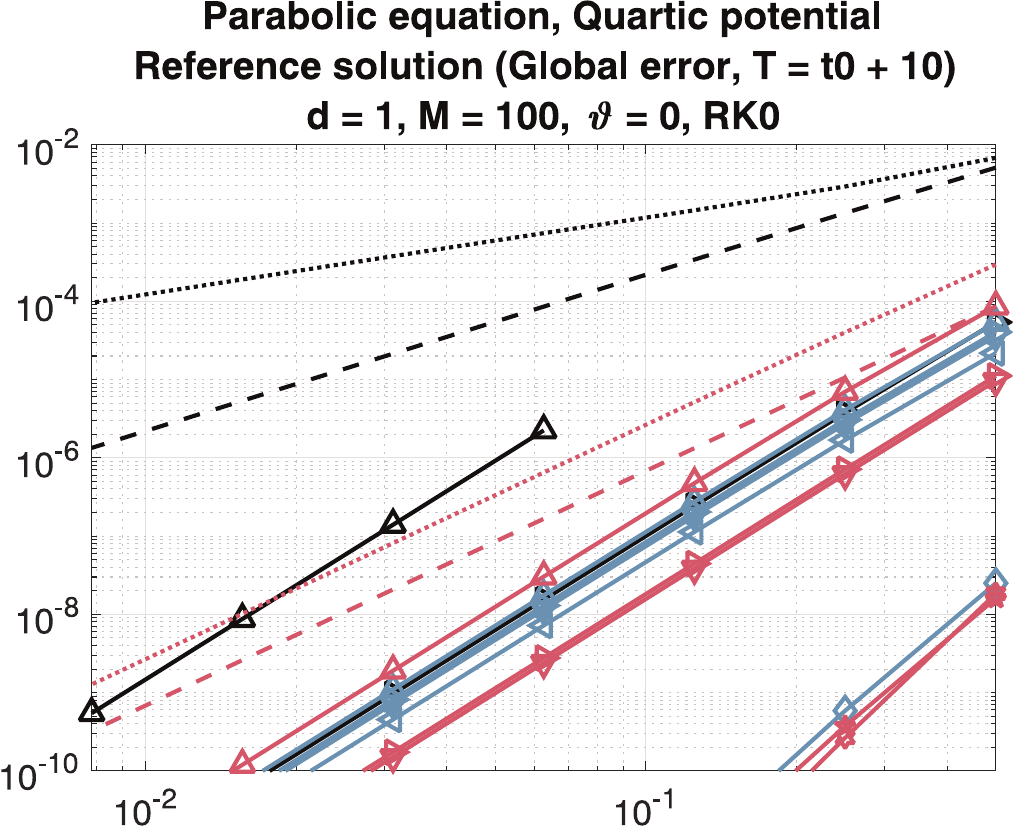} 
\includegraphics[width=0.49\textwidth]{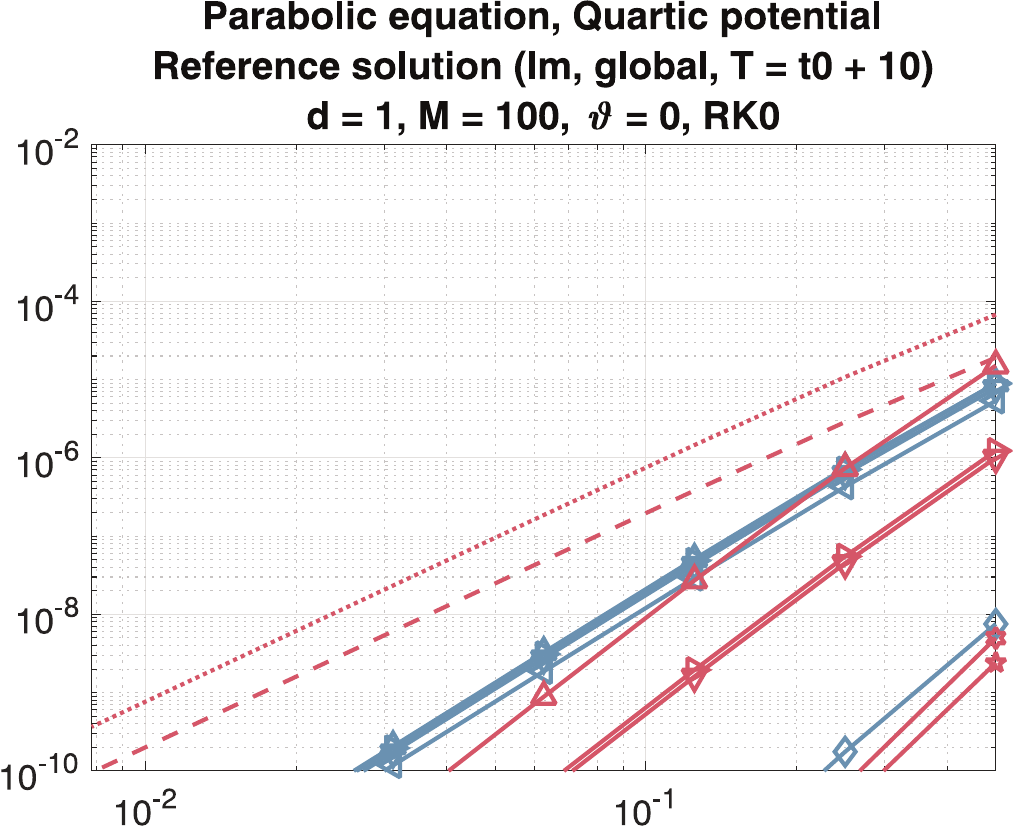} 
\caption{Time integration of the linear parabolic model problem with real-valued solution by real and complex splitting methods, see also~\eqref{eq:ModelProblems} and Figure~\ref{fig:Figure1}.
For the considered quartic potential, a numerical reference solution is computed.
Left: Local and global errors.
Right: Corresponding errors in the imaginary parts.}
\label{fig:Model2_2}
\end{center}
\end{figure}

Altogether, the results displayed in Figures~\ref{fig:Model2_1} and~\ref{fig:Model2_2} show that the nonstiff orders of convergence are retained.
For the errors in the imaginary parts of symmetric-conjugate schemes of even orders, we indeed observe superconvergence, see Table~\ref{tab:Table2} and Section~\ref{sec:Section3}.

\begin{figure}[t!]
\begin{center}
\includegraphics[width=0.47\textwidth]{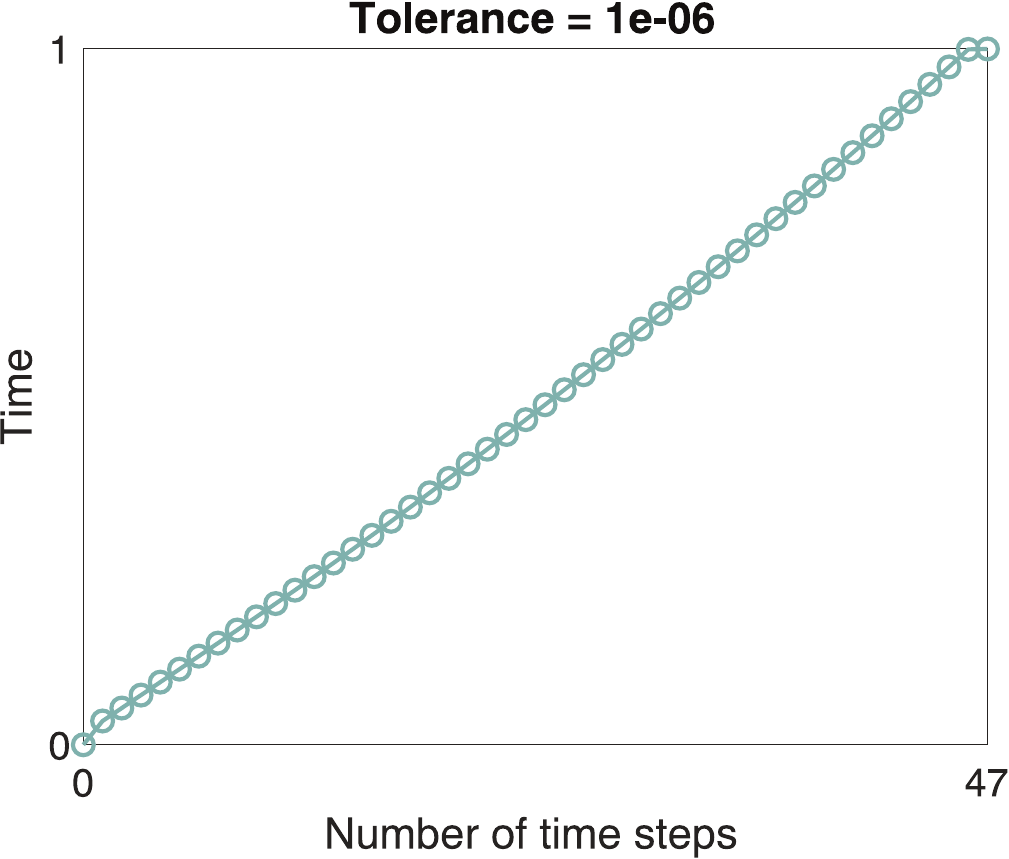} 
\includegraphics[width=0.499\textwidth]{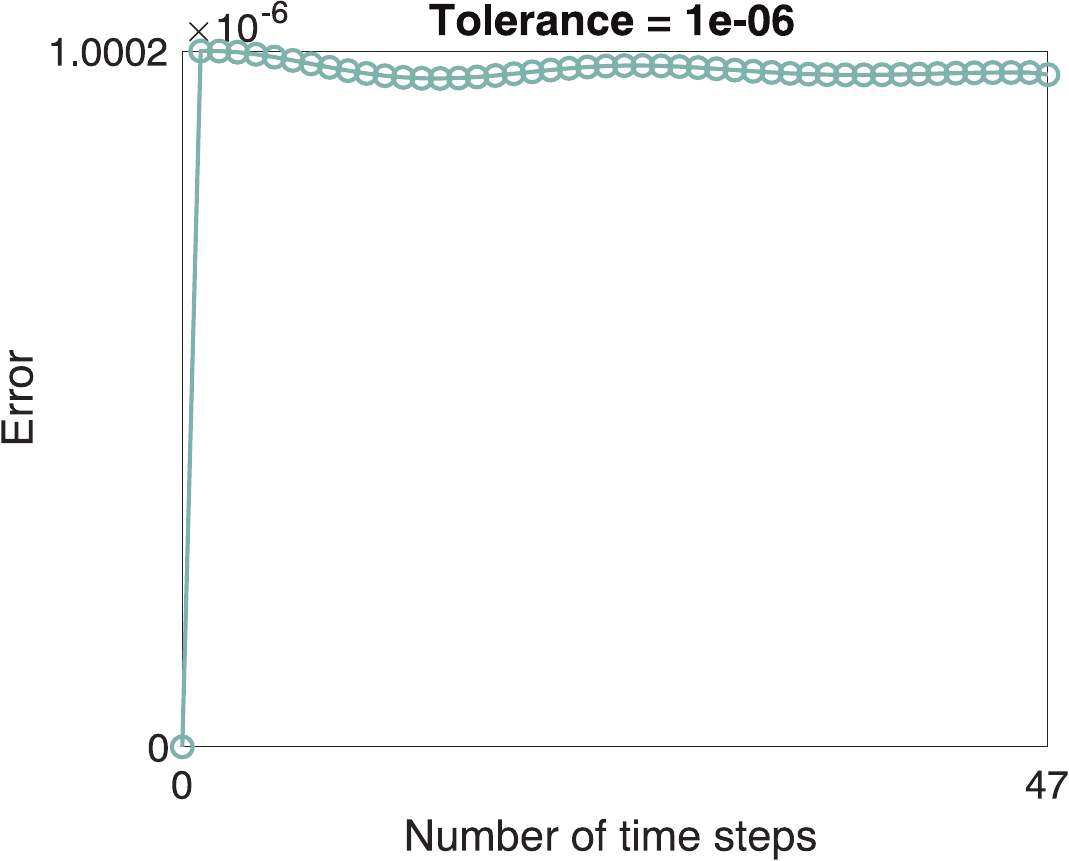} \\[2mm]
\includegraphics[width=0.47\textwidth]{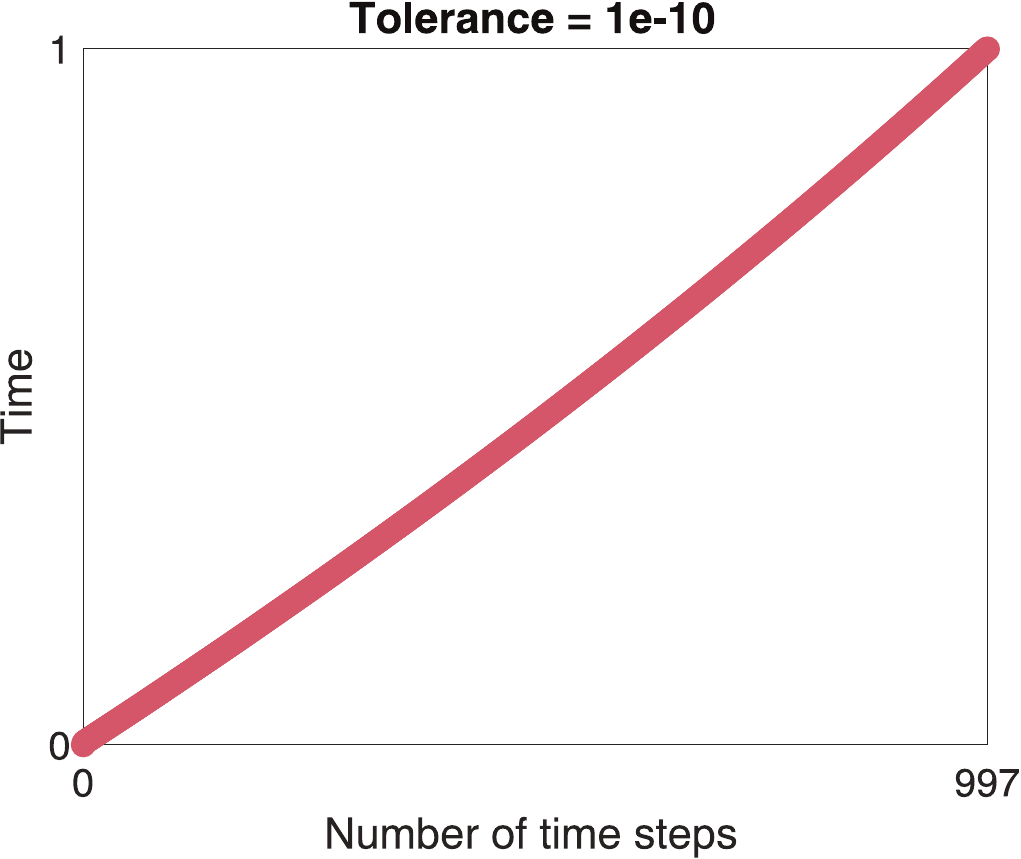} 
\includegraphics[width=0.499\textwidth]{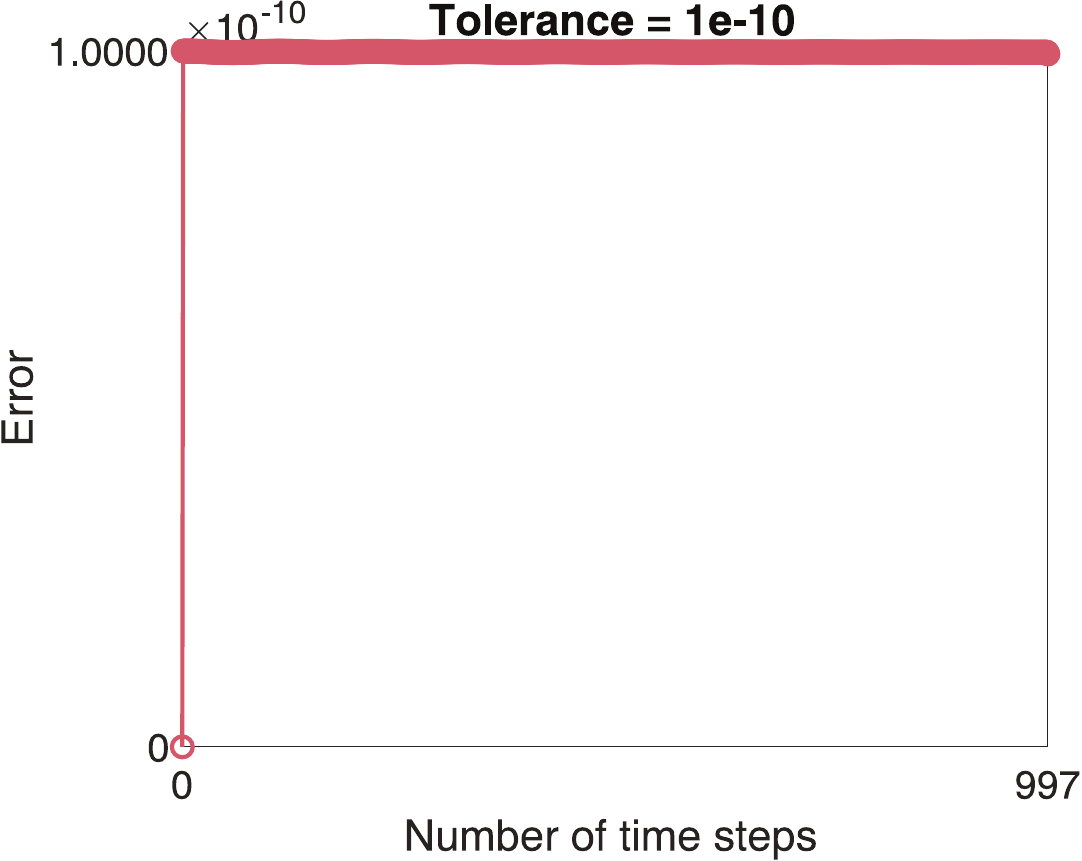} 
\caption{Adaptive time integration of a linear parabolic model problem with known real-valued solution for $t_0 = 0$ and $T = 1$ by a third-order symmetric-conjugate splitting method, see also~\eqref{eq:ModelProblems} and Figure~\ref{fig:Figure1}.
The local error estimation is based on the computation of the imaginary parts with respect to the Euclidean norm and requires negligible additional costs.
The total numbers of time steps 47 and 997 are adjusted in accordance with the prescribed tolerances $10^{-6}$ and $10^{-10}$, respectively.
Left: Sequences of time grid points.
Right: Associated sequences of local errors determined with respect to the exact solution values.}
\label{fig:Model2_3}
\end{center}
\end{figure}

\begin{figure}[t!]
\begin{center}
\includegraphics[width=0.47\textwidth]{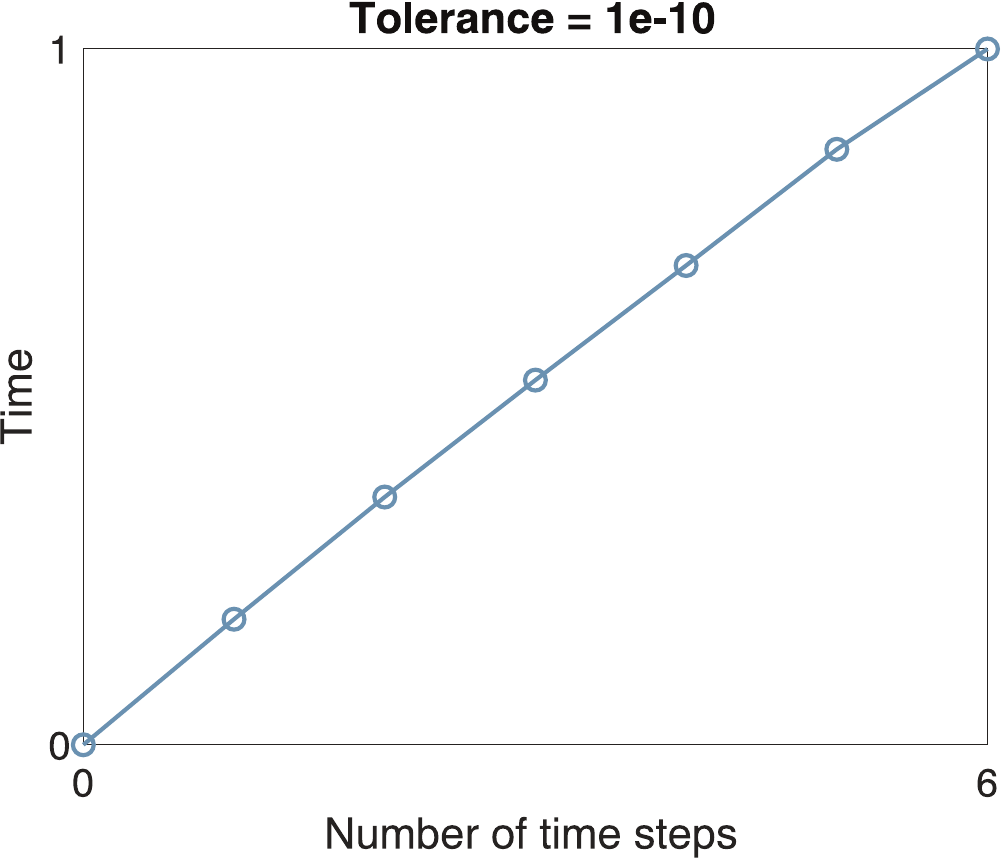} 
\includegraphics[width=0.499\textwidth]{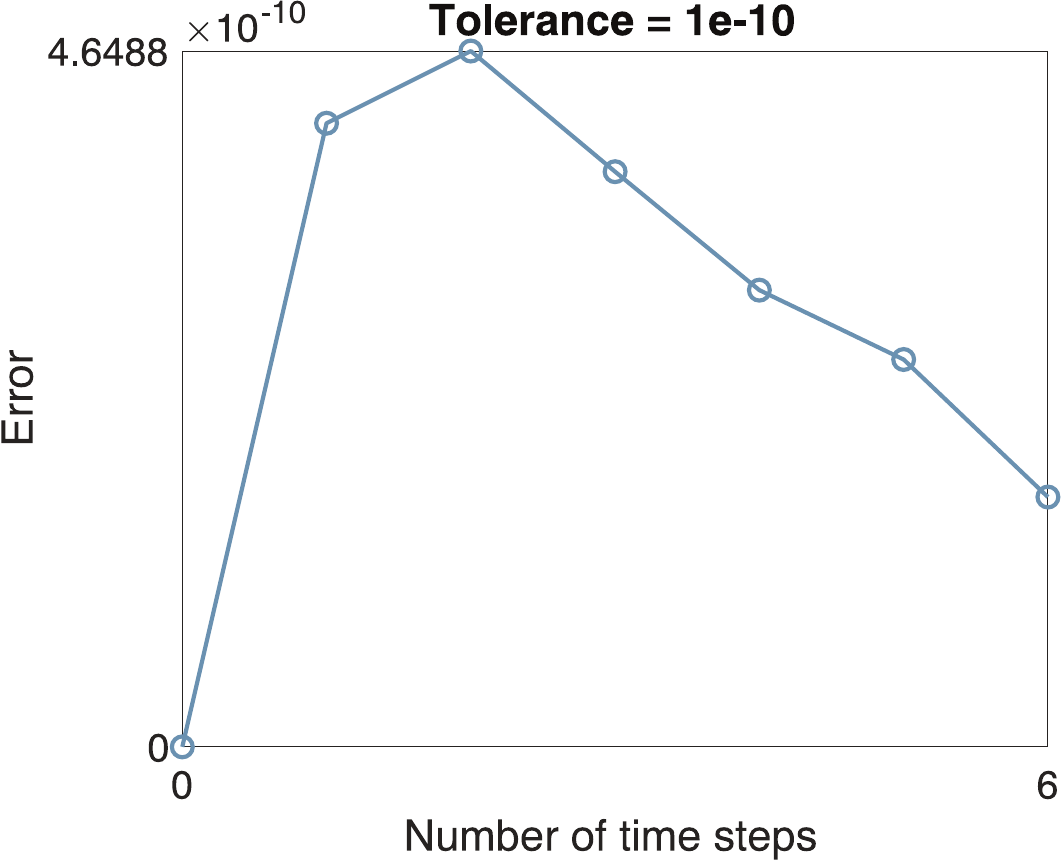} \\[2mm]
\includegraphics[width=0.47\textwidth]{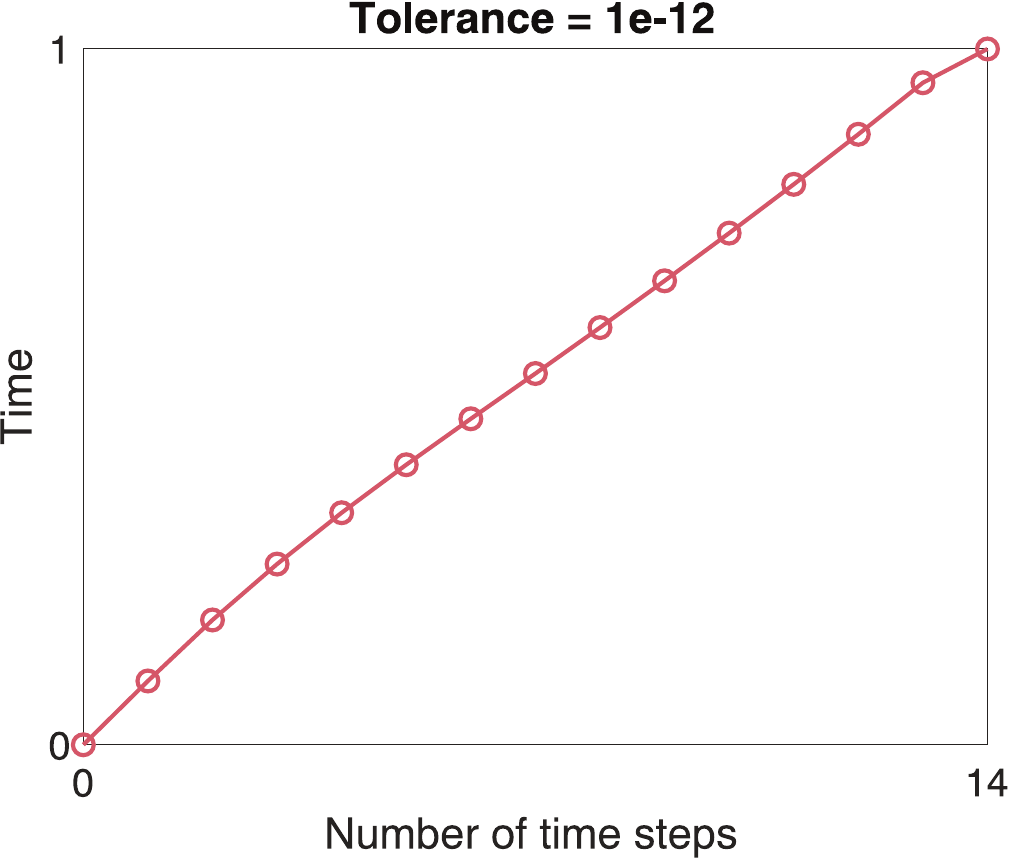} 
\includegraphics[width=0.499\textwidth]{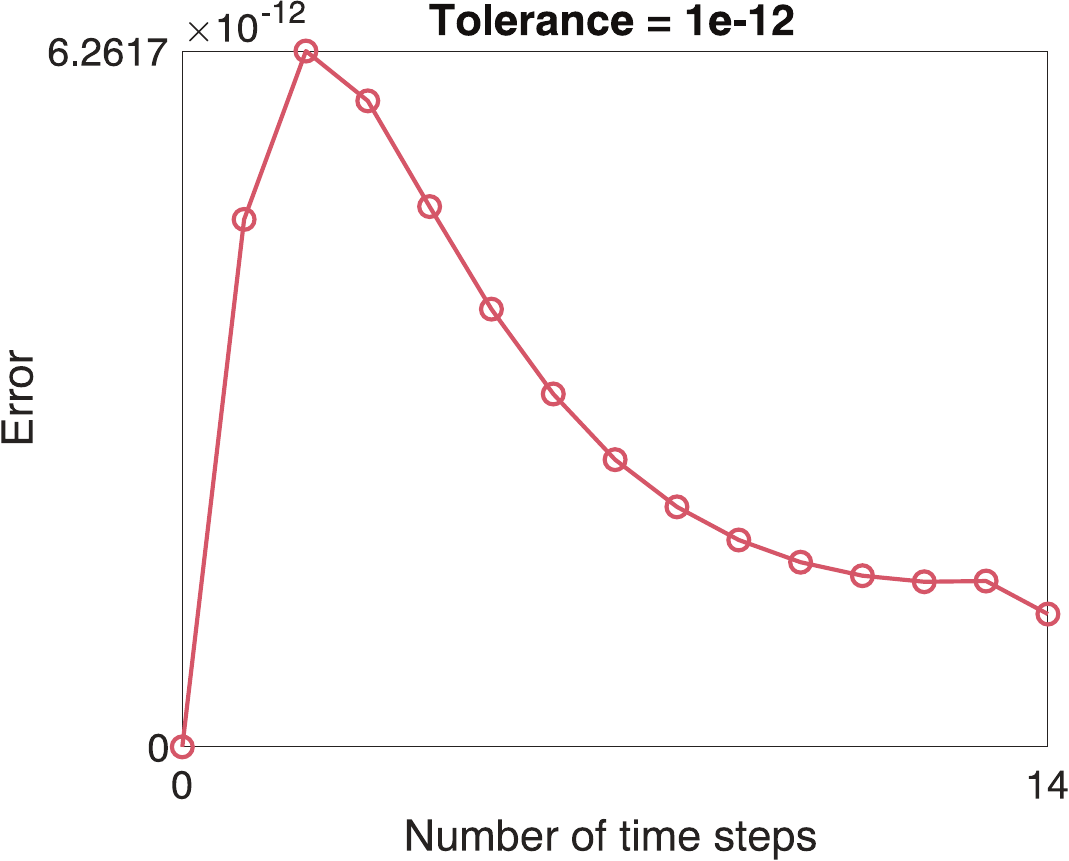} 
\caption{Corresponding results for a sixth-order symmetric-conjugate splitting method and lower tolerances $10^{-10}$ and $10^{-12}$, respectively.}
\label{fig:Model2_4}
\end{center}
\end{figure}

\begin{figure}[t!]
\begin{center}
\includegraphics[width=0.47\textwidth]{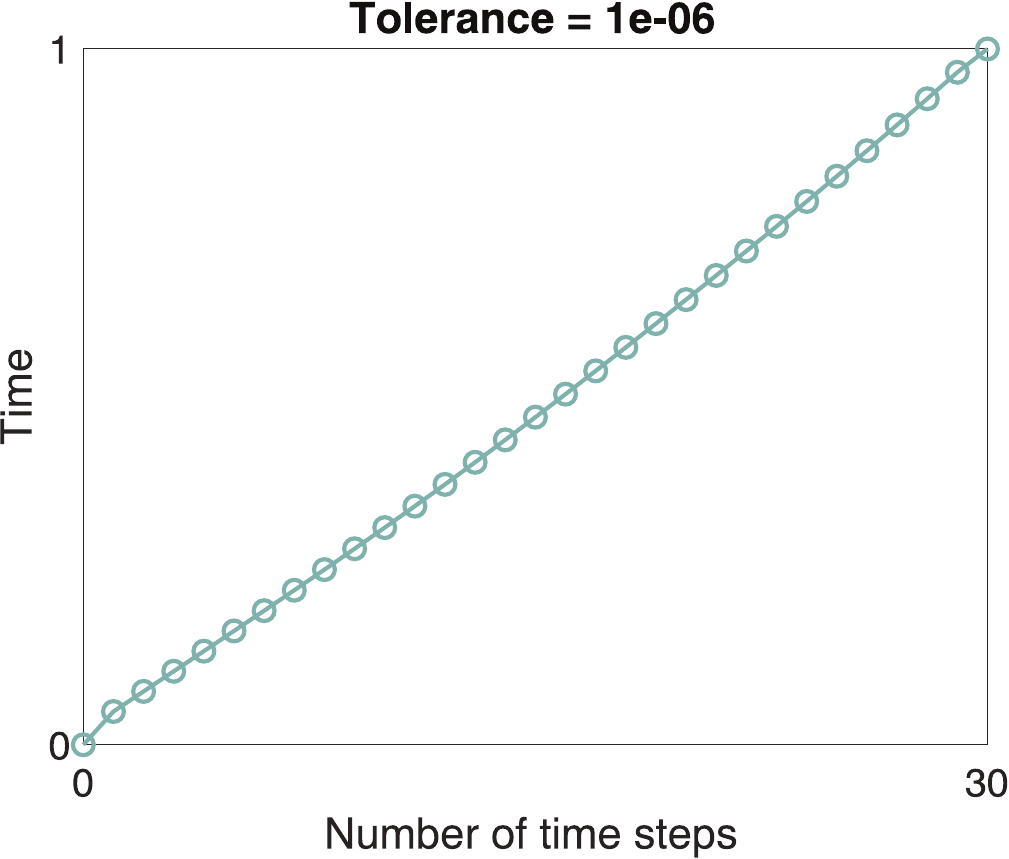} 
\includegraphics[width=0.499\textwidth]{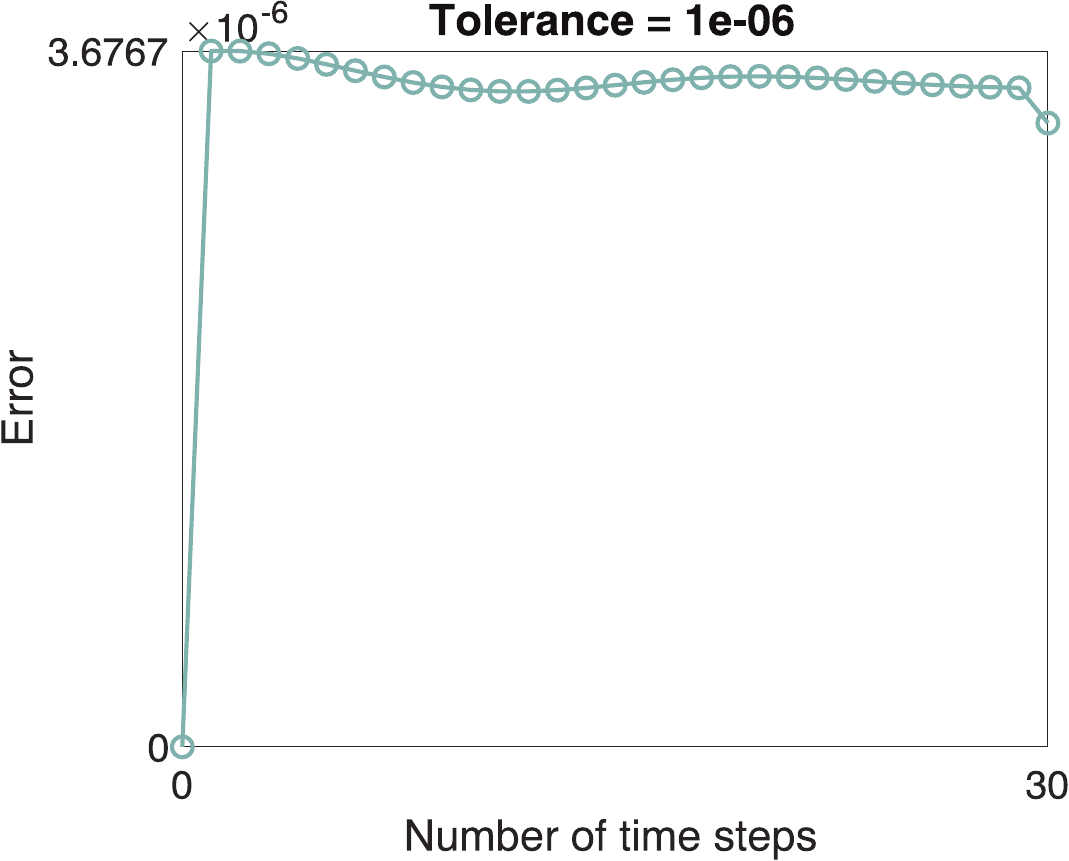} \\[2mm]
\includegraphics[width=0.47\textwidth]{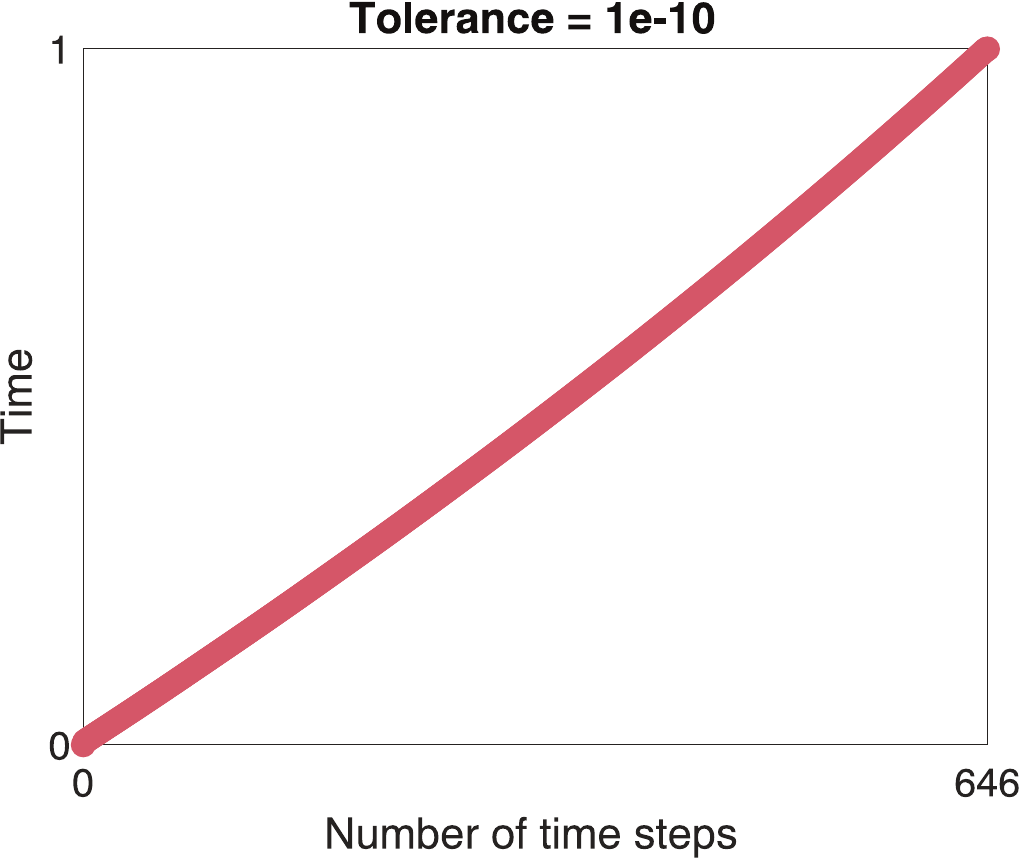} 
\includegraphics[width=0.499\textwidth]{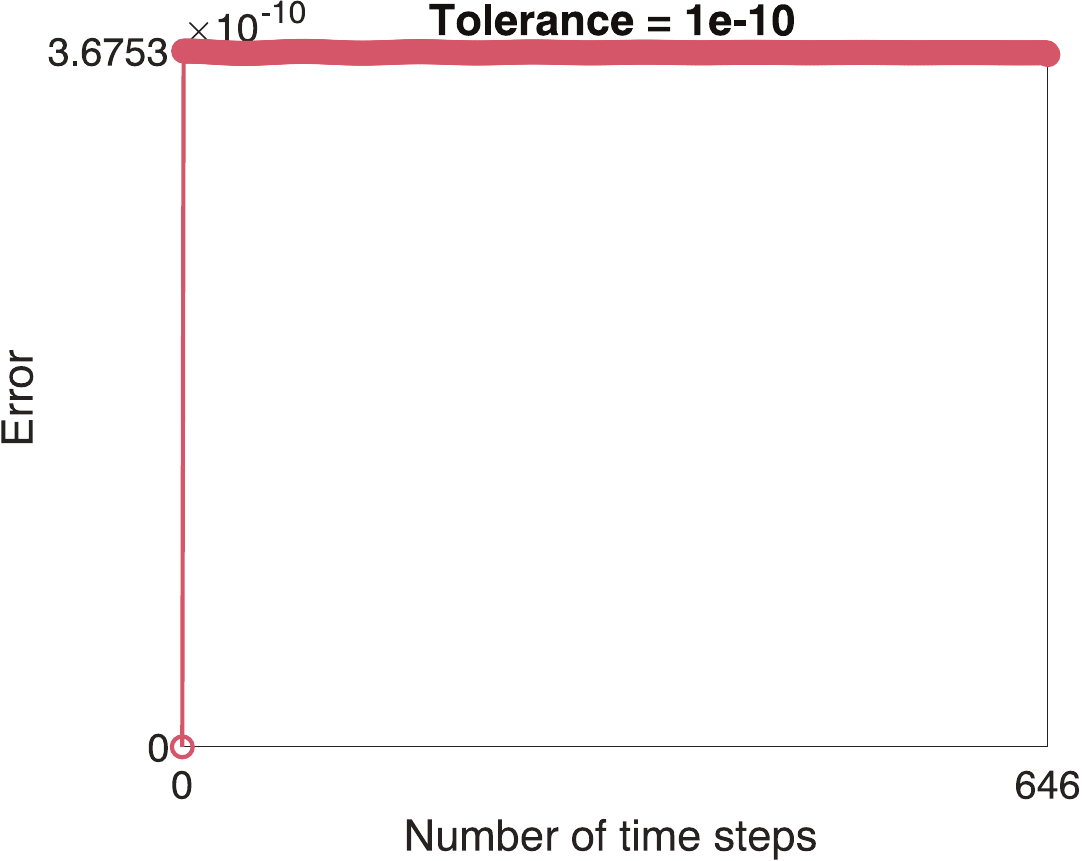} 
\caption{Corresponding results for the third-order symmetric-conjugate splitting method with respect to the maximum norm.}
\label{fig:Model2_5}
\end{center}
\end{figure}

\begin{figure}[t!]
\begin{center}
\includegraphics[width=0.47\textwidth]{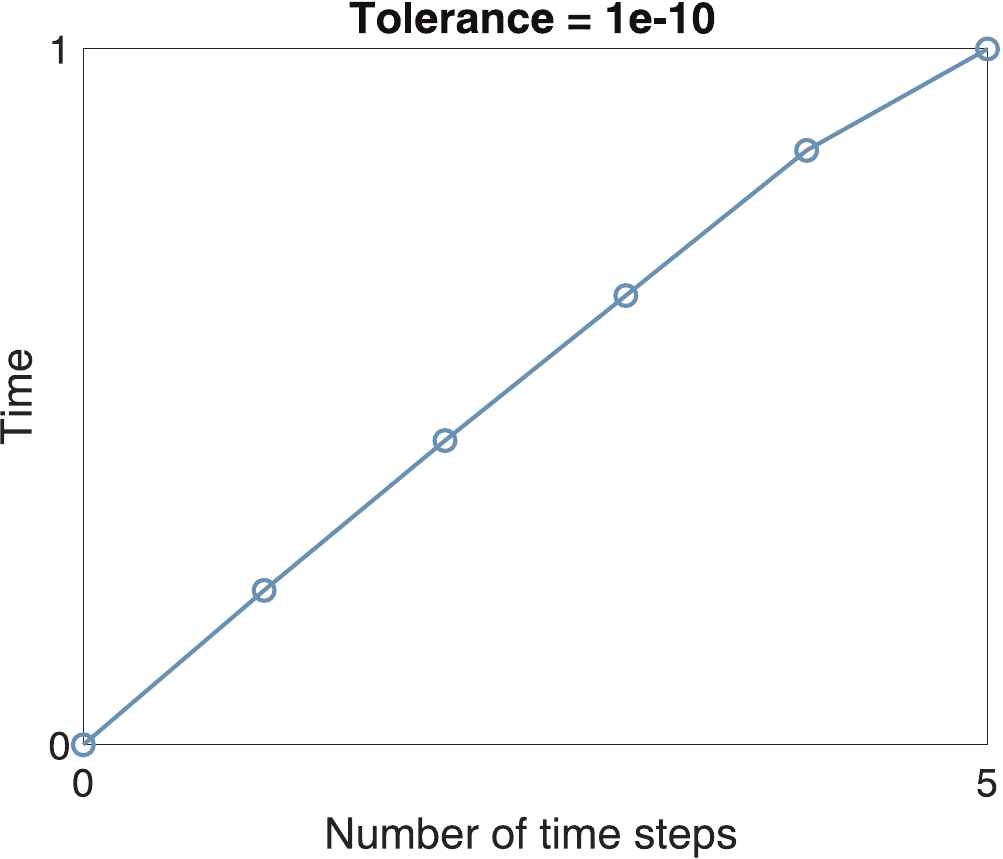} 
\includegraphics[width=0.499\textwidth]{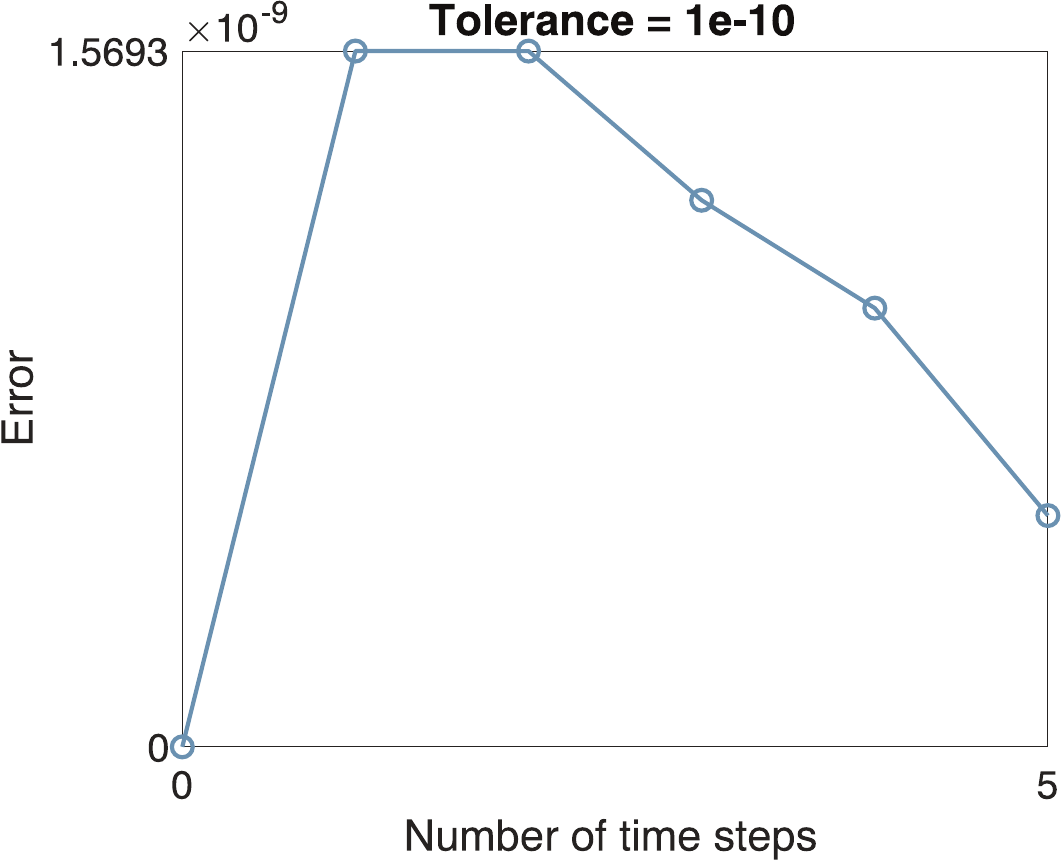} \\[2mm]
\includegraphics[width=0.47\textwidth]{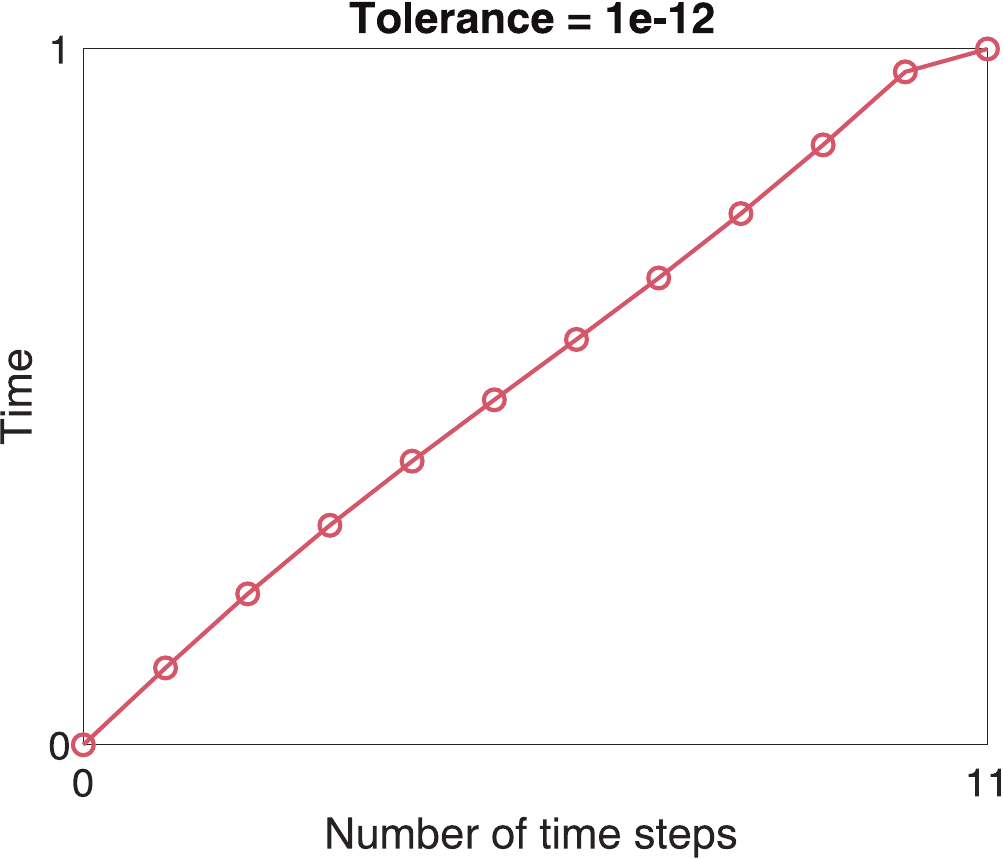} 
\includegraphics[width=0.499\textwidth]{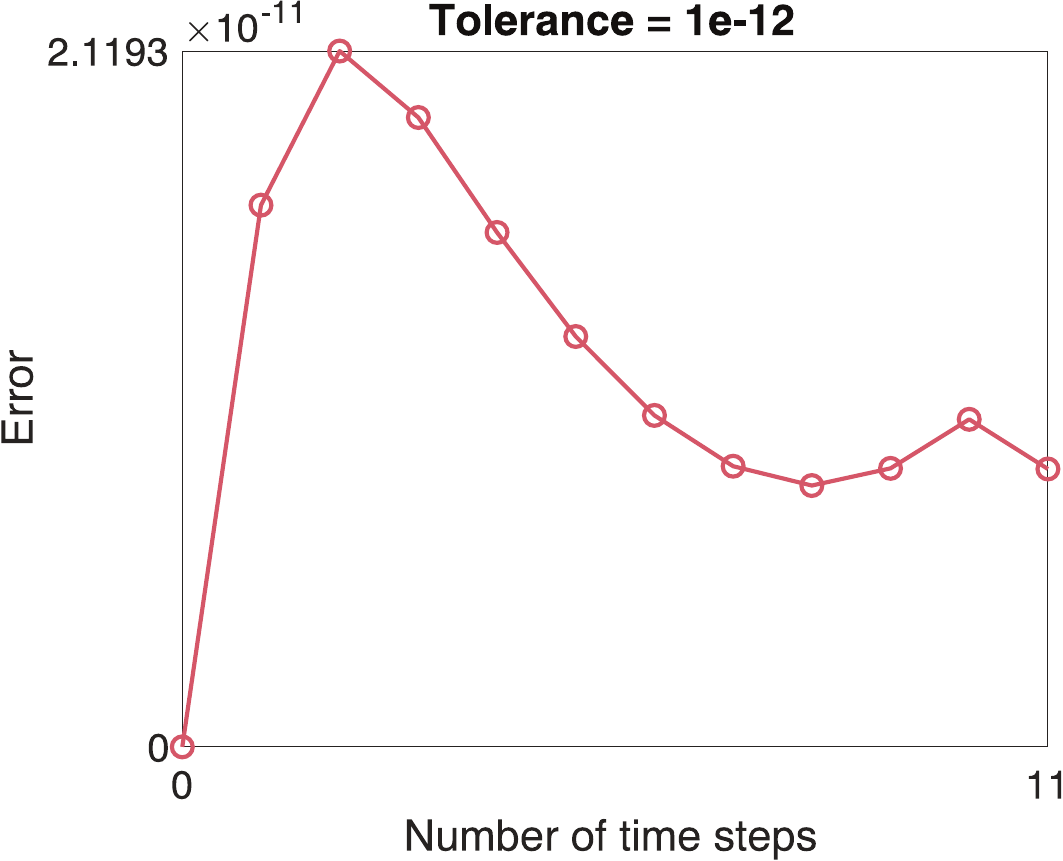} 
\caption{Corresponding results for the sixth-order symmetric-conjugate splitting method with respect to the maximum norm.}
\label{fig:Model2_6}
\end{center}
\end{figure}

\begin{figure}[t!]
\begin{center}
\includegraphics[width=0.8\textwidth]{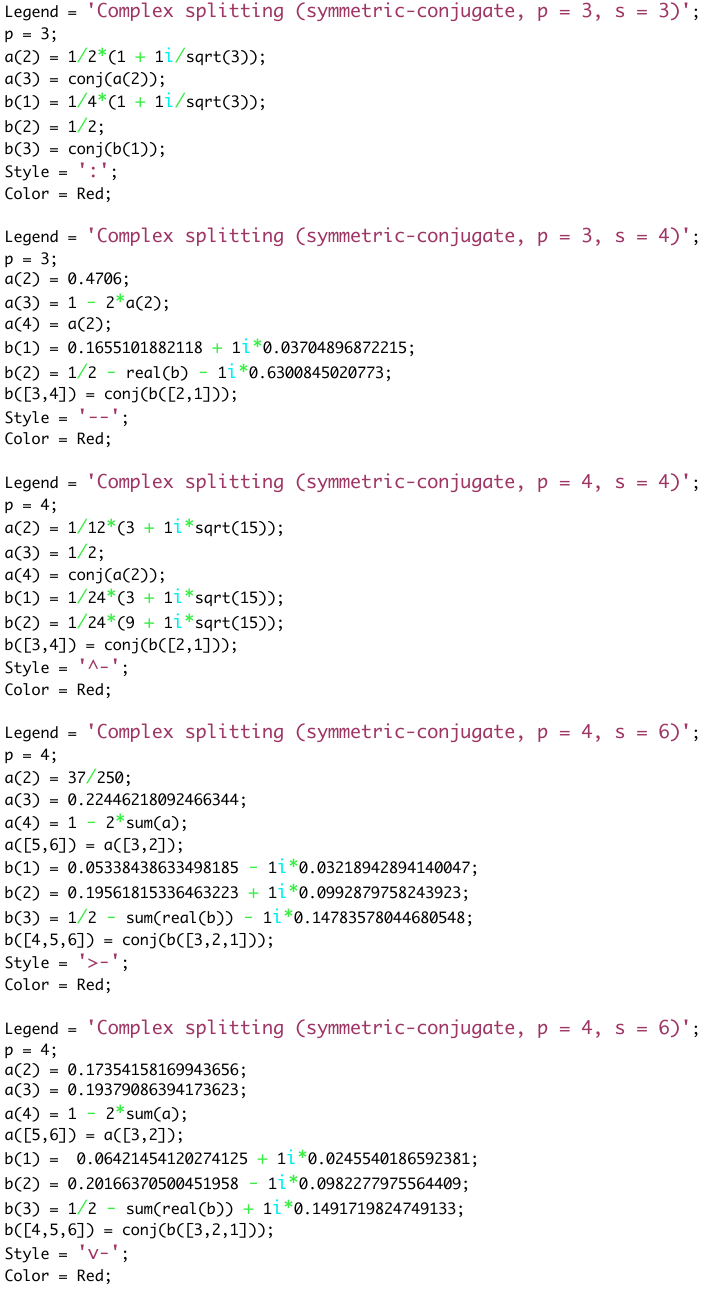} 
\caption{Coefficients of symmetric-conjugate operator splitting methods applied in numerical tests.}
\label{fig:Figure2}
\end{center}
\end{figure}

\begin{figure}[t!]
\begin{center}
\includegraphics[width=0.8\textwidth]{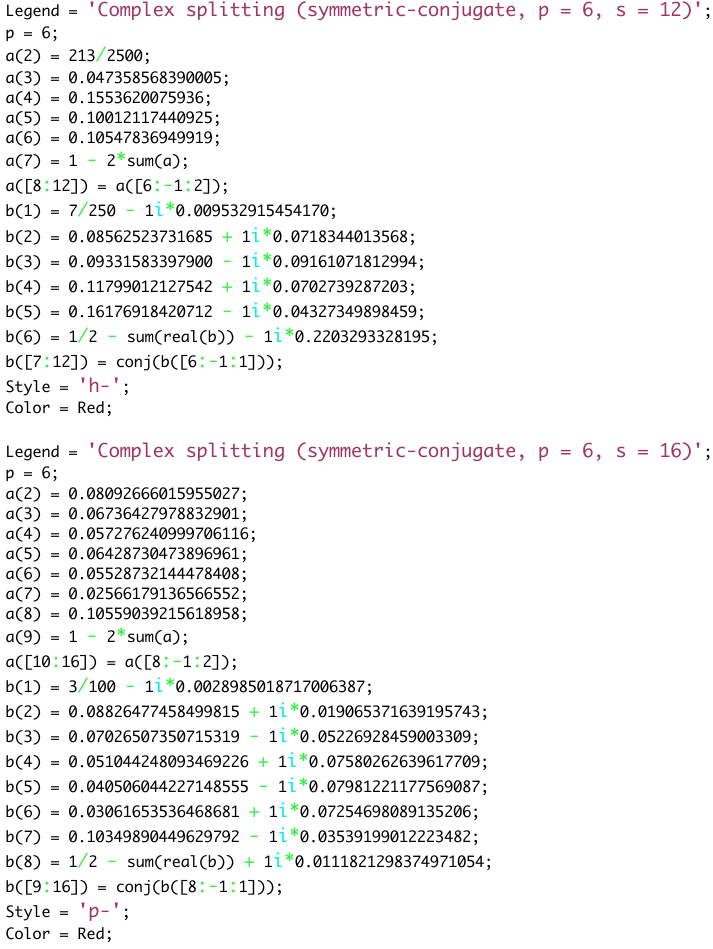} 
\caption{Coefficients of symmetric-conjugate operator splitting methods applied in numerical tests.}
\label{fig:Figure3}
\end{center}
\end{figure}

\MyParagraph{Adaptive time integration}
The design and theoretical analysis of an automatic stepsize control algorithm for the time integration of evolution equations is a complex subject, see for instance~\cite{GustafssonLundhSoederlind1988,HairerNorsetWanner2002,PressTeukolskyVetterlingFlannery2007}.
In order to demonstrate its practical usefulness, we reconsider the linear test equation~\eqref{eq:PDE} for the above specified situation, where the exact solution is available.
A particular benefit of higher-order symmetric-conjugate splitting methods, which distinguishes them from symmetric counterparts, is the possibility to base the local error estimation on the computation of the imaginary parts with negligible additional costs, see Section~\ref{sec:Section3}.
More precisely, for the solution value at the current time~$t_n$, we determine 
\begin{equation*}
\text{local error estimate} = \normbig{q}{\Im\big(u(t_n)\big)}\,, \quad q \in \{2, \infty\}\,.
\end{equation*}
As standard, for a prescribed tolerance, the optimal time stepsize is then adjusted, in essence, through
\begin{equation*}
\tau_{\text{optimal}} = \tau \, \sqrt[p+1]{\frac{\text{tolerance}}{\text{local error estimate}}}\,. 
\end{equation*}
The total number of time steps depends on the choice of the norm, since
\begin{equation*}
\normbig{\infty}{\Im\big(u(t_n)\big)} \leq \normbig{2}{\Im\big(u(t_n)\big)}\,.
\end{equation*}  
Our tests confirm this dependence, yielding smaller time stepsizes for the Euclidean norm ($q = 2$) compared to the maximum norm ($q = \infty$).
The results obtained for third- and sixth-order schemes with respect to the Euclidean norm are displayed in Figures~\ref{fig:Model2_3} and~\ref{fig:Model2_4}.
On the one hand, it is observed that the third-order symmetric-conjugate splitting method performs 47 and 997 time steps to fullfill the requirements that the local errors remain below the prescribed tolerances~$10^{-6}$ and~$10^{-10}$, respectively.
On the other hand, for the sixth-order symmetric-conjugate splitting method and the prescribed tolerances~$10^{-10}$ and~$10^{-12}$, significantly reduced numbers of time steps, namely 6 and 14, are performed.
The corresponding results with respect to the maximum norm are displayed in Figures~\ref{fig:Model2_5} and~\ref{fig:Model2_6}.
Overall, we observe good correlations of the prescribed tolerances for the local errors and the resulting global errors. 
%%%%%%%%%%%%%%%%%%%%%%%%%%%%%%%%%%%%%%%%%%%%%%%%%%%%%%%%%%%%%%%%%%%%%%%%%%%%%%%%%%%%%%%%%%%%%%%%%%%%%%%%%%%%%%%%%%%%%%%%%%%%%%%%%%%%%%%%%%%%%%%%%%%
%%%%%%%%%%%%%%%%%%%%%%%%%%%%%%%%%%%%%%%%%%%%%%%%%%%%%%%%%%%%%%%%%%%%%%%%%%%%%%%%%%%%%%%%%%%%%%%%%%%%%%%%%%%%%%%%%%%%%%%%%%%%%%%%%%%%%%%%%%%%%%%%%%%
\section{Conclusions and future investigations}
%%%%%%%%%%%%%%%%%%%%%%%%%%%%%%%%%%%%%%%%%%%%%%%%%%%%%%%%%%%%%%%%%%%%%%%%%%%%%%%%%%%%%%%%%%%%%%%%%%%%%%%%%%%%%%%%%%%%%%%%%%%%%%%%%%%%%%%%%%%%%%%%%%%
%%%%%%%%%%%%%%%%%%%%%%%%%%%%%%%%%%%%%%%%%%%%%%%%%%%%%%%%%%%%%%%%%%%%%%%%%%%%%%%%%%%%%%%%%%%%%%%%%%%%%%%%%%%%%%%%%%%%%%%%%%%%%%%%%%%%%%%%%%%%%%%%%%%
The present work is dedicated to a comprehensive analysis of symmetric-conjugate operator splitting methods for the time integration of linear evolution equations.
It is seen that the natural approach to incorporate complex coefficients with non-negative real parts permits the design of high-order schemes that remain stable in the context of parabolic problems and thereby overcome the order barriers for standard and modified splitting methods with real coefficients.

Moreover, it is demonstrated that symmetric-conjugate splitting methods are particularly favourable in the numerical integration of nonreversible systems defined by real and symmetric matrices.
The main reasons are that the errors in the imaginary parts and energies remain bounded and hence do not lead to unphysical perturbations. 
Typically, this kind of problems arises in ground and excited state computations for Schr{\"o}dinger equations by the imaginary time propagation, fractal path integrals with applications to many-body theories and statistical physics as well as Monte Carlo simulations of quantum systems.

Future theoretical and numerical investigations will concern extensions to nonlinear evolution equations and a rigorous convergence analysis.
Special attention will be given to complex operator splitting methods applied to complex Ginzburg--Landau equations, since this relevant type of problems interlinks certain characteristics of parabolic as well as Schr{\"o}dinger-type equations.
%%%%%%%%%%%%%%%%%%%%%%%%%%%%%%%%%%%%%%%%%%%%%%%%%%%%%%%%%%%%%%%%%%%%%%%%%%%%%%%%%%%%%%%%%%%%%%%%%%%%%%%%%%%%%%%%%%%%%%%%%%%%%%%%%%%%%%%%%%%%%%%%%%%
%%%%%%%%%%%%%%%%%%%%%%%%%%%%%%%%%%%%%%%%%%%%%%%%%%%%%%%%%%%%%%%%%%%%%%%%%%%%%%%%%%%%%%%%%%%%%%%%%%%%%%%%%%%%%%%%%%%%%%%%%%%%%%%%%%%%%%%%%%%%%%%%%%%
\section*{Acknowledgements}
Part of this work was developed during a research stay at the Wolfgang Pauli Institute Vienna; the authors are grateful to the director Norbert Mauser and the staff
members for their support and hospitality.
This work has been supported by Ministerio de Ciencia e Innovaci{\'o}n (Spain) through projects PID2019-104927GB-C21 and PID2019-104927GB-C22, MCIN/AEI/10.13039/501100011033, ERDF (\emph{A way of making Europe}).
Sergio Blanes and Fernando Casas acknowledge the support of the Conselleria d'Innovaci{\'o}, Universitats, Ci{\`e}ncia i Societat Digital from the Generalitat Valenciana (Spain) through project CIAICO/2021/180.
%%%%%%%%%%%%%%%%%%%%%%%%%%%%%%%%%%%%%%%%%%%%%%%%%%%%%%%%%%%%%%%%%%%%%%%%%%%%%%%%%%%%%%%%%%%%%%%%%%%%%%%%%%%%%%%%%%%%%%%%%%%%%%%%%%%%%%%%%%%%%%%%%%%
%%%%%%%%%%%%%%%%%%%%%%%%%%%%%%%%%%%%%%%%%%%%%%%%%%%%%%%%%%%%%%%%%%%%%%%%%%%%%%%%%%%%%%%%%%%%%%%%%%%%%%%%%%%%%%%%%%%%%%%%%%%%%%%%%%%%%%%%%%%%%%%%%%%

%%%%%%%%%%%%%%%%%%%%%%%%%%%%%%%%%%%%%%%%%%%%%%%%%%%%%%%%%%%%%%%%%%%%%%%%%%%%%%%%%%%%%%%%%%%%%%%%%%%%%%%%%%%%%%%%%%%%%%%%%%%%%%%%%%%%%%%%%%%%%%%%%%%
%%%%%%%%%%%%%%%%%%%%%%%%%%%%%%%%%%%%%%%%%%%%%%%%%%%%%%%%%%%%%%%%%%%%%%%%%%%%%%%%%%%%%%%%%%%%%%%%%%%%%%%%%%%%%%%%%%%%%%%%%%%%%%%%%%%%%%%%%%%%%%%%%%%

\begin{table}[t!]
{\small 
\begin{tabular}{|l|l|}\hline
Lie--Trotter (real, $p = 1$, $s = 1$) & Stability ($a_1 > 0$) \\\hline
Strang (real, $p = 2$, $s = 2$) & Stability ($a_1, a_2 \geq 0$) \\\hline 
Yoshida (real, $p = 4$, $s = 4$) & Instability ($a_3 < 0$) \\\hline\hline
Complex (symmetric, $p = 4$, $s = 4$) & Stability ($\Re \, a_1, \dots, \Re \, a_s \geq 0$) \\\hline
Complex (symmetric, $p = 4$, $s = 5$) & Stability ($a_1, \dots, a_s \geq 0$) \\\hline
Complex (symmetric, $p = 4$, $s = 6$) & Stability ($a_1, \dots, a_s \geq 0$) \\\hline
Complex (symmetric, $p = 6$, $s = 17$) & Stability ($a_1, \dots, a_s \geq 0$) \\\hline\hline
Complex (symmetric-conj., $p = 3$, $s = 3$) & Stability ($\Re \, a_1, \dots, \Re \, a_s \geq 0$)  \\\hline
Complex (symmetric-conj., $p = 3$, $s = 4$) & Stability ($a_1, \dots, a_s \geq 0$)\\\hline
Complex (symmetric-conj., $p = 4$, $s = 4$) & Stability ($\Re \, a_1, \dots, \Re \, a_s \geq 0$) \\\hline
Complex (symmetric-conj., $p = 4$, $s = 6$) & Stability ($a_1, \dots, a_s \geq 0$) \\\hline
Complex (symmetric-conj., $p = 6$, $s = 12$) & Stability ($a_1, \dots, a_s \geq 0$) \\\hline
Complex (symmetric-conj., $p = 6$, $s = 16$) & Stability ($a_1, \dots, a_s \geq 0$) \\\hline
\end{tabular}
}
\bigskip
\caption{Stability properties of real and complex splitting methods in the context of parabolic equations.
Schemes with non-negative coefficients $(a_j)_{j=1}^{s}$ remain stable for Schr{\"o}dinger equations.}
\label{tab:Table1}
\end{table}

\begin{table}[t!]
{\small
\begin{tabular}{|l|c|c|}\hline
Lie--Trotter (real, $p = 1$) & $p_{\, \text{num}} = p = 1$ & --- \\\hline
Strang (real, $p = 2$) & $p_{\, \text{num}} = p = 2$ & --- \\\hline
Yoshida (real, $p = 4$) & $p_{\, \text{num}} = p = 4$ & --- \\\hline\hline
Complex (symmetric, $p = 4$) & $p_{\, \text{num}} = p = 4$ & $p_{\, \text{num}, \Im} = p = 4$ \\\hline
Complex (symmetric, $p = 6$) & $p_{\, \text{num}} = p = 6$ & $p_{\, \text{num}, \Im} = p = 6$ \\\hline\hline
Complex (symmetric-conj., $p = 3$) & $p_{\, \text{num}} = p = 3$ & $p_{\, \text{num}, \Im} = p = 3$ \\\hline
Complex (symmetric-conj., $p = 4$) & $p_{\, \text{num}} = p = 4$ & $p_{\, \text{num}, \Im} = p + 1 = 5$ \\\hline
Complex (symmetric-conj., $p = 6$) & $p_{\, \text{num}} = p = 6$ & $p_{\, \text{num}, \Im} = p + 1 = 7$ \\\hline
\end{tabular}
}
\bigskip
\caption{Application of real and complex splitting methods to parabolic model problems with real-valued solutions.
List of classical orders (first column), numerically observed orders of convergence for solution values (second column), and numerically observed orders for imaginary parts (third column).}
\label{tab:Table2}
\end{table}

%%%%%%%%%%%%%%%%%%%%%%%%%%%%%%%%%%%%%%%%%%%%%%%%%%%%%%%%%%%%%%%%%%%%%%%%%%%%%%%%%%%%%%%%%%%%%%%%%%%%%%%%%%%%%%%%%%%%%%%%%%%%%%%%%%%%%%%%%%%%%%%%%%%

%%%%%%%%%%%%%%%%%%%%%%%%%%%%%%%%%%%%%%%%%%%%%%%%%%%%%%%%%%%%%%%%%%%%%%%%%%%%%%%%%%%%%%%%%%%%%%%%%%%%%%%%%%%%%%%%%%%%%%%%%%%%%%%%%%%%%%%%%%%%%%%%%%%

%%%%%%%%%%%%%%%%%%%%%%%%%%%%%%%%%%%%%%%%%%%%%%%%%%%%%%%%%%%%%%%%%%%%%%%%%%%%%%%%%%%%%%%%%%%%%%%%%%%%%%%%%%%%%%%%%%%%%%%%%%%%%%%%%%%%%%%%%%%%%%%%%%%

%%%%%%%%%%%%%%%%%%%%%%%%%%%%%%%%%%%%%%%%%%%%%%%%%%%%%%%%%%%%%%%%%%%%%%%%%%%%%%%%%%%%%%%%%%%%%%%%%%%%%%%%%%%%%%%%%%%%%%%%%%%%%%%%%%%%%%%%%%%%%%%%%%%
%%%%%%%%%%%%%%%%%%%%%%%%%%%%%%%%%%%%%%%%%%%%%%%%%%%%%%%%%%%%%%%%%%%%%%%%%%%%%%%%%%%%%%%%%%%%%%%%%%%%%%%%%%%%%%%%%%%%%%%%%%%%%%%%%%%%%%%%%%%%%%%%%%%
\end{document}